\documentclass[]{amsart}
\usepackage{amsfonts,amssymb,epsfig}
\usepackage[latin1]{inputenc}
\usepackage[all]{xy}
\usepackage{amsmath}
\usepackage{delarray}
\usepackage{amssymb}
\usepackage{diagrams}
\theoremstyle{definition}

\renewcommand{\epsilon}{\varepsilon}

\newcommand{\cF}{\mathcal F}

\theoremstyle{remark}

\numberwithin{equation}{section}

\newcommand{\cE}{\mathcal E}
\newcommand{\cD}{\mathcal D}
\newcommand{\cG}{\mathcal G}

\newcommand{\cT}{\mathcal T}

\newcommand{\bbN}{\mathbb N}

\newcommand{\R}{{\mathbb R}}

\newcommand{\IE}{\text\textquestiondown}

\newcommand{\cN}{\mathcal N}

\newcommand{\sfD}{\mathsf D}
\newcommand{\sfE}{\mathsf E}
\newcommand{\sfT}{\mathsf T}

\newcommand{\sfS}{\mathsf S}
\newcommand{\sfNc}{\mathsf N_c}
\newcommand{\sfNo}{\mathsf N_o}

\newcommand{\Hom}{\mathsf{Hom}}
\newcommand{\un}{\mathsf{I}}

\newcommand{\nul}{{\scriptscriptstyle^\circ}\!}
\newcommand{\mul}{{\,{\scriptscriptstyle\stackrel{\ast}{\times}}\,}}
\newcommand{\cast}{\square\hskip-6.7pt\ast}
\newcommand{\tast}{\bigtriangleup\hskip-8pt\ast}

\newarrow{ClopTo}{sqhook}---{vee}
\newarrow{OpTo}{hooka}---{vee}
\newarrow{To}----{vee}

\newarrow{iTo}{triangle}---{vee}
\newarrow{ITo}{blacktriangle}---{vee}
\newarrow{sTo}----{triangle}
\newarrow{STo}----{blacktriangle}

\newarrow{iLine}{triangle}----
\newarrow{ILine}{blacktriangle}----
\newarrow{Equal}=====%
\newarrow{DotTo}....{vee}
\newarrow{DotsTo}....{triangle}
\newarrow{DotSTo}....{blacktriangle}
\newarrow{DotiTo}{triangle}...{vee}
\newarrow{DotITo}{blacktriangle}...{vee}

\newarrow{DashTo}{}{dash}{}{dash}{vee}
\newarrow{DashsTo}{}{dash}{}{dash}{triangle}
\newarrow{DashSTo}{}{dash}{}{dash}{blacktriangle}
\newarrow{DashiTo}{triangle}{dash}{}{dash}{vee}
\newarrow{DashITo}{blacktriangle}{dash}{}{dash}{vee}

\newarrow{DashLine}{}{dash}{}{dash}{}
\newarrow{DashiLine}{triangle}{dash}{}{dash}{}
\newarrow{DashILine}{blacktriangle}{dash}{}{dash}{}
\def\rsum{\raisebox{2pt}{$\:\scriptscriptstyle+\:$}}
\def\lsum{\raisebox{2pt}{$\:\scriptscriptstyle+\:$}}
\def\dsum{\vphantom b\scriptscriptstyle+}
\def\usum{\vphantom b\scriptscriptstyle+}

\newarrowfiller{+}\rsum\lsum\dsum\usum
\newarrow{SumSTo}{}{+}{}{+}{blacktriangle}

\def\rrnd{\raisebox{1.5pt}{$\:\scriptscriptstyle{\cdot}\:$}}
\def\lrnd{\raisebox{1.5pt}{$\:\scriptscriptstyle{\cdot}\:$}}
\def\drnd{\vphantom b\scriptscriptstyle{\cdot}}
\def\urnd{\vphantom b\scriptscriptstyle{\cdot}}

\newarrowfiller{rnd}{\rrnd}{\lrnd}{\drnd}{\urnd}
\newarrow{RndSTo}{}{rnd}{}{rnd}{blacktriangle}
\newarrow{RndLine}{}{rnd}{}{rnd}{}

\def\rbul{\raisebox{1.5pt}{$\:\scriptscriptstyle{\bullet}\:$}}
\def\lbul{\raisebox{1.5pt}{$\:\scriptscriptstyle{\bullet}\:$}}
\def\dbul{\vphantom b\scriptscriptstyle{\bullet}}
\def\ubul{\vphantom b\scriptscriptstyle{\bullet}}

\newarrowfiller{bul}{\rbul}{\lbul}{\dbul}{\ubul}
\newarrow{BulSTo}{}{bul}{}{bul}{blacktriangle}
\newarrow{BulLine}{}{bul}{}{bul}{}

\newcommand{\ractTo}[1]{\rTo^{#1}_{\act}}

\newcommand{\rdeqvsTo}[1]{\rdsTo^{#1}_{\eqv}}
\newcommand{\ldeqvsTo}[1]{\ldsTo^{#1}_{\eqv}}
\newcommand{\reqvSTo}[1]{\rSTo^{#1}_{\eqv}}

\newcommand{\rdeqvSTo}[1]{\rdSTo^{#1}_{\eqv}}
\newcommand{\ldeqvSTo}[1]{\ldSTo^{#1}_{\eqv}}

\def\tr{\pitchfork}

\def\eqv{{\fboxsep1pt\fbox{\tiny{$\sim$}}}}
\def\fid{\text{{\fboxsep1pt\framebox{\tiny{\sffamily fid}}}}}
\def\act{\text{{\fboxsep1pt\framebox{\tiny{\sffamily act}}}}}

\def\exa{\text{{\fboxsep1pt\framebox{\tiny{\sffamily exa}}}}}

\def\pb{\text{\raisebox{1.6pt}{{\fboxsep0.7pt\framebox{\tiny{\sffamily pb}}}}}}
\def\po{\text{\raisebox{1.6pt}{{\fboxsep0.7pt\framebox{\tiny{\sffamily po}}}}}}

\def\ful{\text{\raisebox{1.4pt}{{\fboxsep0.5pt\framebox{\tiny{\sffamily ful\phantom{\hskip-3pt{p}}}}}}}}

\def\sepb{\text{\ \tiny{\SEpbk}}}

\newcommand{\SEI}{\begin{diagram}[size=0.5em,tight,abut]
&&\\
{}&&\\
&\rdITo&\\
&&{\cdot}\\
\end{diagram}}

\newcommand{\SES}{\begin{diagram}[size=0.5em,tight,abut]
&&\\
{}&&\\
&\rdSTo&\\
&&{\cdot}\\
\end{diagram}}

\def\et{\text{\textcircled{\scriptsize{\'{e}}}}}

\def\rc{\text{\textcircled{\scriptsize{c}}}}

\newcommand{\gPdS}[2]{\begin{diagram}[size=2em,tight,labelstyle=\scriptstyle,inline]
#1&\pile{\rSTo^{\alpha\ \ }\\ \rSTo_{\beta\ \ }}&#2\\
\end{diagram}}
\newcommand{\gPds}[2]{\begin{diagram}[size=2em,tight,labelstyle=\scriptstyle,inline]
#1&\pile{\rsTo^{\alpha\ \ }\\ \rsTo_{\beta\ \ }}&#2\\
\end{diagram}}

\begin{document}
\title[Generalized atlases]{Lie groupoids as generalized atlases}
\author{Jean PRADINES}
\address{26, rue Alexandre Ducos, F31500 Toulouse, France}
\email{jpradines@wanadoo.fr} 


\date{04/11/03}

\begin{abstract} 
Starting with some motivating examples (classsical atlases for a manifold, space of leaves of a foliation, group orbits), we propose to view a Lie groupoid as a generalized atlas for the ``virtual structure'' of its orbit space, the equivalence between atlases being here the smooth Morita equivalence. This ``structure'' keeps memory of the isotropy groups and of the smoothness as well. To take the smoothness into account, we claim that we can go very far by retaining just a few formal properties of embeddings and surmersions, yielding a very polymorphous unifying theory. We suggest further developments.
\end{abstract}
\maketitle 

\setcounter{tocdepth}{2}
\tableofcontents

\section{Introduction.}
\label{int}The aim of the present lecture is, rather than to present new \emph{results}, to sketch some {\emph{unifying concepts}} and {\emph{general methods}} wished to be in Charles Ehresmann's spirit.

As usual I am expecting that \textit{geometers} will think these sorts of concepts are too general and too abstract for being useful, while \textit{categoricists} will estimate they are too special and too concrete for being interesting. However let us go.

\par\medskip
In the following, I shall be concerned with a certain structure $B$ (basically thought as a manifold) endowed with a certain equivalence relation denoted by  $\sim$ or $R$, and I would like to describe what kind of smoothness or structure is inherited from $B$ by the quotient set $Q=B/R$.
The canonical projection will be denoted by $B\stackrel{q}{\rightarrow}Q$.
The relation $R$ will be identified with its graph, defined by the following pullback
\begin{diagram}[size=1.7em,tight,labelstyle=\scriptstyle,textflow]
        R     &\rTo^{\beta}            &         B   \\
\dTo^{\alpha}         &            {\pb}    & \dTo_{q}       \\
        B     &\rTo_{q}            &         Q   \\
  \end{diagram}%
square, in which $\beta=\text{pr}_1$, $\alpha=\text{pr}_2$. We also denote by $R\stackrel{\tau_R\ }{\rightarrow}B\times B$ the canonical injection, with $\tau_R=(\beta,\alpha)$. 

In case when the structure $B$ is just a topology, the well known answer is given by the so-called quotient or identification topology on $Q$, which owns the good expected universal property in the category \textbf{Top}.
However we notice that, when given other similar data $B'\stackrel{q'}{\rightarrow}Q'$, one has not in general, in spite of a famous error (in Bourbaki's first edition), a homeomorphism between the product 
$Q\times Q'$ and the quotient space of $B\times B'$ by the product of the two equivalence relations, though this is true in two important cases, when $q$ and $q'$ are both open or proper, since $q\times q'$ has the same property.

On the opposite, when $B$ is a manifold, it is well known that there is no such satisfactory answer when staying inside the category $\cD=\textbf{Dif}$ of (smooth maps between) smooth manifolds, i.e. there is no suitable manifold structure for $Q$.

Now for facing this situation there may be two opposite, or better complementary, styles of approaches.

\par\smallskip
The first one consists in ``completing'' $\cD$, i.e. embedding $\cD$ in a larger category $ \widehat{\cD}$ by adding new objects in such a way than $\widehat{\cD}$ has better categorical properties, i.e. has enough limits for allowing to define a good universal quotient. For instance one can wish $\widehat{\cD}$ be a topos.

Various interesting solutions do exist, the study of which is out of our present scope. We just mention, besides Ehresmann's approaches, two dual ways (considered, under various aspects, by several lecturers at the present Conference) of defining generalized smooth structures on $Q$, one (first stressed by Frölicher) consisting in defining the smooth curves, while the second method (emphasized by Souriau with his diffeologies) considers the smooth functions on $Q$. Alain Connes' ``non-commutative'' approach is also related.

\par\smallskip
We follow here an opposite path, avoiding to add too many (necessarily pathological) new objects, and trying to stay within $\cD$. We do not attempt to define a generalized smooth structure (in the set-theoretical sense) on the most general quotients, and limit ourself to objects which are sufficiently close to manifolds in the sense that they can be described by means of equivalence classes of some simple types of diagrams in $\cD$ ; we do not try to introduce the limits of such diagrams in the categorical sense.

\par\medskip
\emph{Indeed we think that the classical categorical concept of limit involves in general a certain loss of the information encapsulated in the concept of a suitable equivalence class of diagrams, but we shall not attempt to develop more formally such a general concept here, though we think it a very promising way, being content with illustrating this point of view by the important special case sketched presently.}

\section{Some motivating examples.}
\label{Motex}
Before going to abstract general definitions, I start by giving some elementary examples (to be made more precise later) of the kinds of objects I have in mind.

\subsection{Regular equivalences.}
\label{Regeq}
The ideal situation is of course that of the so-called regular equivalences. This means that there exists on $Q$ a (necessarily unique) manifold structure such that $B\rSTo^{q\ \ } Q$ is a surmersion ($=$ surjective submersion). (Here we start anticipating some pieces of notation for arrows to be systematized later within a more general setting).

\emph{Godement's theorem} gives a characterization of those equivalences by properties of the graph $R$ summarized by the following notations :

\begin{center}
\fbox{{\gPdS{R}{B}}\ \ and\ \ $R\rITo ^{\tau_R\ }B\times B$}\,
\end{center}
where again the black triangle head for an arrow stands for ``surmersion'', while the black triangle tail means ``embedding'' (in the sense of Bourbaki), or ``proper embedding'' when dealing with Hausdorff manifolds.

\par\smallskip
These conditions express that $R$, regarded (in a seemingly pedantic way) as a subgroupoid of the (banal) groupoid $B\times B$\,, is indeed a smooth (or Lie) groupoid in the sense introduced by Ehresmann, embedded in $B\times B$, and the manifold $Q$ may be viewed as the \textit{orbit space} of this groupoid.

A Lie groupoid $R$ satisfying the framed conditions will be called a \emph{principal} or \emph{Godement} groupoid.

We shall see in the next example why it is convenient to consider
$$\gPdS{R}{B}\rSTo^{q\ \ } Q$$ as a ``\emph{generalized \emph{(non \'{e}tale)} atlas}'' for $Q$.

If we have another ``atlas'' $\gPdS{R'}{B'}\rSTo^{q'\ \ } Q$ of the same manifold $Q$, we can take the fibred product of $q$ and $q'$, and we get a commutative diagram 
\begin{diagram}[size=2.8em,tight,labelstyle=\scriptstyle,p=0.7em,textflow,midshaft]%
 & &S & & \\
 &\ldeqvSTo{f} &\dDashSTo^{\beta_S}\dDashSTo_{\alpha_S} &\rdeqvSTo{f'}& \\
R & & E& & R' \\
\dSTo^{\beta_R}\dSTo_{\alpha_R}&\ldSTo^{\quad p} &  &\rdSTo^{p'\quad} &\dSTo^{\beta_{R'}}\dSTo_{\alpha_{R'}}\\
B & &\dDashSTo_r & & B' \\
 &\rdSTo_q & &\ldSTo_{q'} & \\ 
 & & Q & & \\
\end{diagram}
expressing the ``\emph{compatibility}'' of these ``atlases'', which means that they define the same (manifold) structure on $Q$.

More precisely the graph $S$ (which in turn may be viewed as a groupoid) can be obtained in the following way by means of the commutative cube below,
\begin{diagram}[size=2em,loose,labelstyle=\scriptstyle,textflow]
{S}& &\reqvSTo{f}& &R & & \\
   &\rdeqvSTo{f'}\rdDashSTo(6,2)^h_{\eqv} & & &\dILine &\rdeqvSTo g & \\
\dDashITo^{\tau_S} & &R' &\reqvSTo{g'} &\HonV & &\, Q \, \\
   & &\dITo^{\tau_{R'}}& &\dTo^{\tau_R} & & \\
{E\times E}&\hLine^{p\times p}&\VonH&\rSTo&{B\times B} & &\dITo^{\tau_Q} \\
   &\rdSTo^{p'\times p'}\rdDashSTo(6,2)^{r\times r} & & & &\rdSTo^{q\times q}& \\
   & &B'\times B'& &\rSTo^{q'\times q'}& &Q\times Q \\
\end{diagram}%
 the bottom face of which is the pullback of $q\times q$ and $q'\times q'$, and which is constructed step by step by pulling back along the vertical arrows, starting with $\tau_Q$. The last one is just the diagonal of $Q$, and may be considered as the anchor map of the ``\emph{null}'' groupoid $Q$ (consisting of just units). The upper face is then also a pullback, as well as all the six faces, and also the vertical diagonal square with three dashed edges. The $\sim$ symbols emphasize (very special instances of) ``surmersive equivalences'' between Lie groupoids. (One can observe on this diagram the general property of ``parallel transfer by pulling back'' for the embeddings and surmersions).

Thus we see that the ``compatibility'' of the two (generalized) atlases $$\gPdS{R}{B}\rSTo^{q\ \ } Q \mathrm{\ and\ } \gPdS{R'}{B'}\rSTo^{q'\ \ } Q$$ for $Q$ is expressed by the existence of a common ``refinement'' (pictured above with dashed arrows) : $\gPdS{S}{E}\rSTo^{r\ \ } Q$.

One might prove directly a converse, which indeed follows from more general considerations.

\subsection{Classical atlases.}
\label{Classat}
An important special case of the previous one (which explains the terminology) will be given by the following diagrammatic description of atlases and covers of a manifold $Q$.
\par\medskip
Let $(\varphi_i:V_i\rightarrow U_i)_{(i\in I)}$ be an atlas of the manifold $Q$, where $(V_i)_{(i\in I)}$ is a cover of $Q$ and the codomains $U_i$ 's of the charts $\varphi_i$ 's are open sets in some model space (which may be $\R^n$ or a Banach space).

Let $V_{ij}$ be $V_i\cap V_j$ and $U_{ij}$ be the image of $V_{ij}$ in $U_i$ by the restriction of $\varphi_i$.

Set $V=\coprod_{i\in I}V_i$, with its canonical projection $r:V\rightarrow Q$ (whose datum is equivalent to the datum of the covering), $U=\coprod_{i\in I}U_i$ (a trivial manifold), $R=\coprod_{(i,j)\in I\times J}U_{ij}$, and $S=\coprod_{(i,j)\in I\times J}V_{ij}$. The charts $\varphi_i$ 's define a bijection $\varphi:V\rightarrow U$ as well as a bijection $\phi:S\rightarrow R$.

Note that $S$, together with its canonical projection onto $Q$, defines the intersection covering, while, with its two canonical projections onto $V$, it can be viewed also as the graph of the equivalence relation associated to the surjection $r$.

Using the bijections $\varphi$, $\phi$, we have analogous considerations for $R$ and $U$, but moreover the latter are (trivial) manifolds, and the equivalence is regular, so that we recover a (very) special instance of the situation in the first example. Here the projections $\alpha_R$, $\beta_R$ are not only surmersions but moreover 
\emph{\'{e}tale} maps (of a special type, which might be called trivial) ; 
here they will be pictured by arrows of type $\rsTo$. More precisely their restrictions to the components of the coproduct $R$ define homeomorphisms onto the open sets $U_{ij}$ 's, and the datum of the smooth groupoid $R$ with base $U$ is precisely equivalent to the datum of the \emph{pseudogroup of changes of charts}.

The situation is summed up by the following diagram (with $q=r\circ\varphi^{-1}$),
\begin{diagram}[size=2em,tight,p=0.6em,labelstyle=\scriptstyle,textflow]%
S&\rDotTo^{\phi}_{\approx}&R\\
\dDotTo^{\beta_S}\dDotTo_{\alpha_S}& &\dsTo^{\beta_R}\dsTo_{\alpha_R}\\
V&\rDotTo^{\varphi}_{\approx}&U\\
\dDotTo^r& &\dsTo^q\\
Q&\rEqual&Q\\
\end{diagram}%
 which describes the \emph{generalized atlas $\gPds{R}{U}\rsTo^{q\ \ } Q$
associated to a classical atlas} (whence the terminology). The dotted arrows in the diagram are to remind that the left column lies in \textbf{Set}, while the right column lies in \textbf{Dif}.

\par\smallskip
Now if we define a refinement of the previous atlas, denoted by $U\rsTo^{q\ \ }Q$ for brevity, as an atlas $W\rsTo^{r\ \ }Q$ such that $r$ admits a surjective
\footnote{This surjectivity is not implied by the usual definition of a refinement of a covering, but one can always impose it by the following slight modification of the definition, which changes nothing for the common use made of it : one demands in addition that a refinement of a covering contains this latter covering, which can always be achieved by taking their union.}
 factorization $W\rsTo^{p\ \ } U$, it is easy to see that the compatibility of two atlases $U\rsTo^{q\ \ } Q$ and $U'\rsTo^{q'\ \ } Q$ may be expressed by the existence of a common refinement $W\rsTo^{r\ \ } Q$, and we get a special case of the notion of compatibility introduced in the previous subsection, where the general surmersions are replaced by (very special) surjective \'{e}tale maps.

The compatibility diagram now reads as below, 
\begin{diagram}[size=2.8em,tight,labelstyle=\scriptstyle,p=0.7em,textflow]%
 & &Z & & \\
 &\ldeqvsTo{f} &\dsTo^{\beta_Z}\dsTo_{\alpha_Z} &\rdeqvsTo{f'}& \\
R & & W& & R' \\
\dsTo^{\beta_R}\dsTo_{\alpha_R}&\ldsTo^{\quad p} &  &\rdsTo^{p'\quad} &\dsTo^{\beta_{R'}}\dsTo_{\alpha_{R'}}\\
U & &\dsTo_r & & U' \\
 &\rdsTo_q & &\ldsTo_{q'} & \\ 
 & & Q & & \\
\end{diagram}
and this explains the terminology of the previous section.

\subsection{Group and groupoid actions.}
\label{GA}
An action of a Lie group $G$ on a manifold $E$, defined by the (smooth) map : $G\times E\stackrel{\beta}{\rightarrow}E$, $(g,x)\mapsto g\cdot x=\beta(g,x)$, can be described by the graph of the map $\beta$. Modifying the order in the products, this graph defines an embedding ${H=G\times E}\rITo^{\iota}_{} {G\times(E\times E)}$, $(g,x)\mapsto(g,(g\cdot x,x))$. Regarding $G\times(E\times E)$ as a (Lie) groupoid with base $E$ (product of the group $G$ by the ``banal'' groupoid $E\times E$), the associativity property of the action law may be expressed by the fact that $H$ is an (embedded) subgroupoid of $G\times(E\times E)$. Composing $\iota$ with ${\text{pr}_1}$ yields a (smooth) functor $H\stackrel{f}{\rightarrow}G$ which owns the property that the commutative diagram generated by the source projections is a pullback.

This construction extends for the action of a (smooth) groupoid \gPdS{G}{B} acting on a manifold over $B$ : $E\rSTo^{p\ }B$, replacing the product $G\times(E\times E)$ by the fibred product of the anchor map $\tau_G$ and $p\times p$, and one gets a pullback square :
\begin{diagram}[size=2.2em,tight,labelstyle=\scriptstyle]
        H     &\rSTo^f            &         G   \\
\dSTo^{\alpha_H}         &            {\pb}    & \dSTo_{\alpha_G}       \\
        E     &\rSTo^{f^{(0)}}            &         B   \\
  \end{diagram} %
where $f^{(0)}=p$ (induced by the functor $f$ on the bases $E=H^{(0)}$\,, $B=G^{(0)}$). The action law may be recovered by $\beta=\beta_H$, using the isomorphism of $H$ with the fibred product of $G$ and $E$ over $B$, which results from the pullback property.

Note that the pullback property remains meaningful even when $E\stackrel{p}{\rightarrow}B$ is not a surmersion, since $\alpha_G$ still is, (though the fibred product of $\tau_G$ by $p\times p$ may then fail to exist), and $\beta_H$ still defines an action of $G$ on $p$.

For that reason we call a functor $H\ractTo{f} G$ owning the previous pullback property an \emph{actor}, which is emphasized by the framed label of the arrow. Such functors received various unfortunate names in the categorical literature, among which ``discrete opfibrations'' and ``foncteurs d'hypermorphisme'' (Ehresmann), and, better, ``star-bijective'' (Ronnie Brown), but note that the present concept encapsulates a smoothness information, included in the pullback property, and not only the purely set-theoretic or algebraic conditions (see below for a more general setting).

There is an equivalence of categories between the category of equivariant maps between action laws and the category admitting the actors as objects and commutative squares of functors as arrows.

Note than in the literature $H$ is currently called the \emph{action groupoid}, but it is only the whole datum of the actor $f:H\rightarrow G$ which fully describes the action law, whence our terminology.

\par\smallskip
Here we let $Q$ be the (set-theoretic) quotient of the manifold $E$ by the action of the Lie group(oid) $G$. By the previous construction it appears too as the orbit space of the Lie groupoid $H$, so we have again for $Q=E/G=E/H$ a generalized atlas $\gPdS{H}{E}\rDotSTo Q$, the dotted arrow meaning here that we have now just a set-theoretic surjection (we have here to go out of {\textbf{Dif}}, since $Q$ is no more a manifold).

\subsection{Foliations on $B$.}
\label{Fol}
Here we need a more restrictive notion for our surmersions, called \emph{retroconnected} (it is in a certain sense precisely the opposite of being \'{e}tale, which might be called as well ``\emph{retrodiscrete}''), and, in the present subsection, unlike the previous one, a notation such as $E\rsTo^p B$ will indicate that the surmersion p is retroconnected. This means that the inverse image of any $x\in B$ is connected, or, equivalently, that the inverse image of any connected subset of $B$ is connected.

In the following, when we have to use simultaneously \'{e}tale and retroconnected surmersions, we shall distinguish them by means of circled labels :
\begin{center}
$A\rsTo_{\et}B$\ \ \ or\ \ \ $A\rsTo_{\rc}B$\ .
\end{center}

\par\smallskip
There is an obvious notion of foliation \emph{induced} by pulling back a foliation along such a retroconnected surmersion, and the induced foliation keeps the same set-theoretical space of leaves.

Then the notion of {\sffamily F}-equivalence (in the sense introduced by P. Molino in the 70's) between two foliations $(B,\cF)$, $(B,\cF')$ can be expressed by the existence of a commutative diagram as below, 
\begin{diagram}[size=2em,tight,labelstyle=\scriptstyle,textflow]
      &           &   {(E,\cG)}    &          &        \\
      &\ldsTo^p  &           &\rdsTo^{p'} &        \\
{(B,\cF)}&           &   \dDotsTo^r    &          & {(B,\cF')}  \\
      &\rdDotsTo^{q}  &           &\ldDotsTo^{q'} &        \\
      &           &   Q    &          &        \\
\end{diagram}
which means that $\cF$ and $\cF'$ induce the same foliation $\cG$ on the manifold $E$.

Now we note that again $Q=B/\cF$ may be viewed as the orbit space of a Lie groupoid (here with connected source fibres), to know the Ehresmann \emph{holonomy groupoid} 
\begin{diagram}[size=1.8em,tight,labelstyle=\scriptstyle,inline]
H&\pile{\rsTo^{\alpha\ \ }\\ \rsTo_{\beta\ \ }}&B\
\end{diagram} (sometimes renamed much later as the graph of the foliation), and Molino equivalence may alternatively be described by the commutative diagram below 
\begin{diagram}[size=3.3em,tight,labelstyle=\scriptstyle,p=0.5em,textflow]%
 & &K & & \\
 &\ldeqvsTo{f} &\dsTo^{\beta_K}\dsTo_{\alpha_K} &\rdeqvsTo{f'}& \\
H & & E& & H' \\
\dsTo^{\beta_H}\dsTo_{\alpha_H}&\ldsTo^{f^{(0)}} &  &\rdsTo^{f'^{(0)}} &\dsTo^{\beta_{H'}}\dsTo_{\alpha_{H'}}\\
B & &\dDotsTo_r & & B' \\
 &\rdDotsTo_q & &\ldDotsTo_{q'} & \\ 
 & & Q & & \\
\end{diagram}
(with all maps retroconnected), where $f^{(0)}=p$, $f'^{(0)}=p'$, and the symbols $\sim$ mean (surmersive retroconnected) equivalences of Lie groupoids, in a sense to be made more precise later.

As we shall see also later, such a diagram defines a (smooth in a very precise sense) Morita equivalence between $H$ and $H'$.

\medskip
\textbf{Remark.} Any surmersion $A\rSTo^f Q$ admits of an essentially unique\footnote{i.e. up to isomorphisms.} factorization    $$A\rsTo^{q'}_{\rc}Q'\rsTo^e_{\et}Q$$
with $q'$ retroconnected and $e$ \'{e}tale.

Though this might be proved directly, it is better to apply the general theory  of Lie groupoids, in the special case of principal (or Godement) groupoids (see \ref{Regeq} above):

if $R$ is the graph of the regular equivalence on $A$ defined by $q$, its neutral (or $\alpha$-connected) component $R^c$ is an open (possibly non closed\footnote{As in [B], we have to deal with possibly non Hausdorff manifolds.}, but automatically invariant) subgroupoid of $R$, hence it is still a Godement groupoid, with now connected $\alpha$-fibres, and it defines the regular foliation admitting $Q'$\,\footnote{Possibly non Hausdorff.} for its space of leaves. Then $Q'\rTo Q$ is a surmersion with discrete fibres, i.e. \'{e}tale, and the graph of the equivalence on $Q'$ thus defined is isomorphic to the two-sided quotient groupoid $R//R^c$ (cf. [P2]).

\section{Diptychs.}
I am now enough motivated for introducing more dogmatically some general abstract definitions modelling and unifying the previous situations, as well as myriads of others.

\subsection{Definition of diptychs.}
\label{Defdip}
The notions presented in this section have a much wider range than it would be strictly necessary for the sequel, if one wants to stay in \textbf{Dif}, but give to it a much wider scope, even when aiming only at applications in \textbf{Dif}, as illustrated by some of the previous examples.

We introduced them a long time ago, in [P1], and think they deserve being better known and used.
\par\smallskip
In the presentation of the examples of the previous section, we emphasized the role played by embeddings/surmersions (these are special mono/epimorphisms of \textbf{Dif}, but not the most general ones, which would indeed be pathological), with possibly some more restrictive conditions added.

Our claim is that an incredible amount of various constructions can be performed without using the specificity of these conditions, but just a few very simple and apparently mild stability properties (of categorical nature) fulfilled in a surprisingly wide range of situations encountered by the ``working mathematicians''. The power of these properties comes from their conjunction.

\emph{Then the leading idea (illustrated beforehand in the previous section) will be to describe the set-theoretical constructions by means of diagrams, emphasizing injections/surjections, and then rereading these diagrams in the category involved, using the distinguished given mono/epi's.}

\subsubsection{Diptychs data.}
\label{Dipda}
 A ``\emph{diptych}'' $\mathsf{D}=(\cD;\cD_i,\cD_s)$ (which may be sometimes denoted loosely by $\cD$ alone) is defined by the following \emph{data} :
\par\medskip

$\bullet$ $\cD$ is a category which comes equipped with \emph{finite} non void\footnote{See below.} \emph{products}.

The subgroupoid of invertible arrows (called \emph{isomorphisms}) is denoted by $\cD_\ast$.

$\bullet$ $\cD_i$/$\cD_s$ is a subcategory of $\cD$, the arrows of which are mono/epi-morphisms (by axiom (iii) below), called \emph{good mono/epi}'s and denoted generally by arrows with a triangular tail/head such as $\riTo$/$\rsTo$, or $\rITo$/$\rSTo$, and so on (here / is written loosely for resp.).
\par\smallskip
The arrows belonging to $\cD_r=\cD_i\cD_s$ (i.e. composed of a good mono and a good epi) are called \emph{regular}. \emph{In general $\cD_r$ will not be a subcategory}. When $\cD_r=\cD$, the diptych may be called \emph{regular}.

\subsubsection{Diptychs axioms.}
\label{dipax}

These data have to satisfy the following \textit{axioms} (which look nearly self dual, but not fully, and that has to be noticed\footnote{It is said that the same thing happened at the very instant of the big bang, with analogous consequences.}) :
\begin{enumerate}
\item[(i)] $\cD_i\cap\cD_s=\cD_\ast$ ;
\item[(ii)] $\cD_i$ and $\cD_s$ are \emph{stable by products} ;
\item[(iii)] (a)/(b) the arrows of $\cD_i$/$\cD_s$ are \emph{monos/strict}\footnote{This means that they are coequalizers [McL]. See below for an alternative formulation.} \emph{epis} ;
\item[(iv)] (\emph{``strong/weak source/range -stability''} of $\cD_i$/$\cD_s$) :
         \begin{enumerate}
         \item[(a)]\footnote{We just mention that it may be sometimes useful to work with only the ``weak source-stability'' condition, dual of (b). It is then possible to define a suitable \emph{full subcategory} of $\cD$ in which the strong axiom is satisfied. The objects of this subcategory own (in particular) the property that their diagonal maps are in $\cD_i$, which is formally a Hausdorff (or separation) type property. These objects may be called \emph{$i$-scattered}.} $(h=gf\in\cD_i)\Rightarrow(f\in\cD_i)$ ;
         \item[(b)]$((h=gf\in\cD_s)\text{\ \emph{and}\ }(f\in\cD_s))\Rightarrow(g\in\cD_s)$ ;
         \end{enumerate}
\item[(v)](``\emph{transversality}'', denoted by $\cD_s\pitchfork\cD_i$) :
         \begin{enumerate}
         \item[(a)](``\emph{parallel transfer}'') :\\given $A\rSTo^s B$ and $B'\rITo^i B$ (which means : $s\in\cD_s$, $i\in\cD_i$), there exists a pullback
\begin{diagram}[size=2.7em,tight,labelstyle=\scriptstyle,textflow]
        {\IE A'?}     &\rDashITo^{\ \ \IE i'?}            &         A   \\
\dDashSTo^{\IE s'?}         &            {\IE\pb?}    & \dSTo^{}_s       \\
        B'     &\rITo^{}_i            &         B   \\
  \end{diagram}
with moreover $s'\in\cD_s$, $i'\in\cD_i$ (the question marks frame the objects or properties, as well as the dashed arrows, which appear in the conclusions, as consequences of the data).

          \item [(b)](conversely : ``\emph{descent}'', or ``\emph{reverse transfer}'') :\\ given a pullback square as below 
\begin{diagram}[size=2em,labelstyle=\scriptstyle,textflow]
        A' &&    \rITo^{i'} &           &         A   \\
\dSTo^{s'}  &&            {\pb} &   & \dSTo_s       \\
 {\ \ \ \ \  B'\ \ \ \IE}&\rDashILine&{?}   &\rTo_{i\ \ }            &         B   \\
  \end{diagram}
(with $i'\in\cD_i$, $s,\!s'\in\cD_s$, $i\in\cD$), one has $i\in\cD_i$ (conclusion pictured by the question marks around the dashed triangular tail).
            \end{enumerate}
\end{enumerate}

\subsubsection{Full subdiptychs.}
\label{fulsub}
\textbf{Remark.}
Let be given a subclass $\mathcal C$ of the class of objects of $\cD$, which is \emph{stable by products} and let $\cD'$, $\cD'_i$, $\cD'_s$ be the full subcategories of $\cD$, $\cD_i$, $\cD_s$ thus generated.

In order they define a new dyptich $\sfD'$ (full subdiptych), one of the following conditions is sufficient :
\begin{enumerate}
	\item given a good epi $A\rSTo B$, the condition $B\in\mathcal C$ implies $A\in\mathcal C$ ;
	\item given a good mono $A\rITo B$, the condition $B\in\mathcal C$ implies $A\in\mathcal C$ ;
\end{enumerate}

\subsection{Some variants for data and axioms.}
\label{dipvar}
\subsubsection{Prediptychs.}
\label{predip}

When dealing with diagrams, we shall need a weaker notion :
\textbf{Definition.} A \emph{prediptych} is a triple $\mathsf{T}=( \cT; \cT_i, \cT_s)$ where one just demands $ \cT_i$,$ \cT_s$ to be subcategories of $ \cT$\,, such that :
\begin{center}
$ \cT_i\cap \cT_s= \cT_\ast$\ (isomorphisms).
\end{center}
Most of the prediptychs we shall consider are \emph{regular}, i.e. $\cT=\cT_r=\cT_i\cT_s$.
\subsubsection{Alternatives for axioms.}
\label{altax}
For any arrow $B\rTo^f B'$, we have the ``\emph{graph factorization}'' : $B\rITo^i B\times B'\rTo^{\text{pr}_2} B'$,
 with $i=(1_B,f)\in \cD_i$ by axioms (iv) (a) and (i) (as a section of $\text{pr}_1$), whence it is readily deduced that one has $\cD_s\pitchfork\cD$ ; this means :\\
\bigskip
$\bullet$ (v)(a) \emph{remains valid when omitting $\cD_i$ both in assumptions and conclusions}.
\par\medskip
We have also $\cD_s\pitchfork\cD_s$. In particular the pullback square generated by two good epis $p$, $q$, has its four edges in $\cD_s$. Such pullback squares will be called \emph{perfect squares}. The case $p=q:B\rSTo Q$ covers the general situation embracing various examples above.
It can be shown, using composition of pullback squares, that the axiom (iii) (b), in presence of the other ones, may be rephrased in the two following equivalent ways (using the notations of (\ref{Regeq})) :
\begin{enumerate}
\item[$\bullet$] (iii) (b') \emph{any good epi $q:B\rSTo Q$ is the coequalizer of the pair $\gPdS{R}{B}$, where $R$ is the fibred product of $q$ by $q$} ;
\item[$\bullet$] (iii) (b") \emph{any perfect square is a push out too} (this last property is very important and remarkable).
\end{enumerate}

\subsubsection{Terminal object.}
\label{TO}
The existence of the void product (i.e. of a \emph{terminal object}), is not always required, in view of important examples ; when it does exist, it will be denoted by a plain dot $\bullet$ (though its support has not to be a singleton).

Though, in many examples, not only there exists a terminal object, but moreover the canonical arrows $A\rightarrow \bullet$ are in $\cD_s$, it may be useful however not to require this property in general. Then those objects owning this property will be called \emph{$s$-condensed} (see examples below).

If $A\rSTo B$ is in $\cD_s$, then $A$ is $s$-condensed iff $B$ is.

Observe that, if $A$ is $s$-condensed, then, for any object $Z$, the canonical projection $\text{pr}_2:A\times Z\rSTo Z$ will be in $\cD_s$ (which is always true whenever all objects are required to be $s$-condensed).

If $\cD'$, $\cD'_i$, $\cD'_s$ are the \emph{full subcategories generated by the $s$-condensed objects} it can be checked that they define a \emph{(full) subdiptych} $\mathsf{D'}=(\cD';\cD'_i,\cD'_s$).

\subsection{A few examples of basic large diptychs.}
\label{dipex}
In fact most of the categories used by the ``working mathematicians'' own one or several natural diptych structures, and checking the axioms may sometimes be a more or less substantial (not always so well known) and often non trivial part of their theory, which is in this way encapsulated in the (powerful) statement that one gets a diptych structure.\footnote{It might be also advisable to look for a way of adapting the axioms without weakening their power in order to include some noteworthy exceptions such as measurable spaces or Riemannian or Poisson manifolds.}

This is all the more remarkable since a general category (with finite products) bears no canonical non trivial (i.e. with $\cD_s\neq\cD_\ast$) diptych structure, the crucial point being that in general the product of two epimorphisms fails to be an epimorphism.

\subsubsection{Sets.}
 The category $\cE=$ \textbf{Set} of (applications between) sets owns a canonical diptych structure $\mathsf E=\mathsf{Set}$, which is regular, by taking for $\cE_i$/$\cE_s$ the subcategories of injections/surjections (here these are exactly all the mono/epi -morphisms).

The same is true for the dual category (exchanging injections and surjections), but the dual diptych $\mathsf E^\ast$ is not isomorphic to $\mathsf E$\,.\footnote{Owing to the lack of symmetry for the axioms, one cannot in general define the dual of a diptych.}

\subsubsection{Two general examples.}
\label{Abtop}

 There are two remarkable and important cases when one gets a (canonical) diptych structure, which is moreover \emph{regular}, by taking as good monos/epis \emph{all} the mono/epi -morphisms, to know : the \emph{abelian categories} and the \emph{toposes}\footnote{See for instance, for a good part, but not the whole, of the properties involved in this statement, the textbooks by Mac Lane and Peter Johnstone.}. Of course these two very general examples embrace in turn a huge lot of special cases in Algebra and Topology.

As a consequence \emph{all the constructions we carry out by using diptychs are working for general toposes}, but the converse is false, since the most interesting and useful diptychs \emph{are not toposes}.

\subsubsection{Topological spaces.}
\label{top}
 In \textbf{Top} (resp. \textbf{Haus}\footnote{The full subcategory of Hausdorff spaces.}), we can take as good monos the (resp. proper\footnote{If we had taken these ones for \textbf{Top}, we would have been in the situation alluded to in footnote 8 and the Hausdorff spaces might be constructed as the i-scattered objects.}) topological embeddings, and as good epis the surjective \emph{open} maps. All objects are then $s$-condensed. This diptych is \emph{not regular}.

These canonical diptychs will be denoted by {\sffamily Top} and {\sffamily Haus}.

\par\smallskip
We may alternatively take as good epis in \textbf{Top} the \'{e}tale/retroconnected (or, in \textbf{Haus}, proper) surjective maps. Then the \emph{$s$-condensed objects} are the discrete/connected (compact) spaces. This is illustrated by examples above.

\subsubsection{Banach spaces.}
\label{ban}

 In \textbf{Ban}, the category of (continuous linear maps between) Banach spaces, one can take as good monos/epis the left/right invertible arrows. All objects are $s$-condensed. This diptych is not regular (save for the full subdiptych of finite dimensional spaces). It will be denoted by {\sffamily Ban}.

\subsubsection{Manifolds.}
\label{man}

In \textbf{Dif}, the category of (smooth maps between) possibly Banach and possibly non Hausdorff manifolds (in the sense of Bourbaki [B]), the basic diptych structure, denoted by {\sffamily Dif}, is defined by taking as good monos the embeddings and as good epis the surmersions.

But ones gets a very large number of very useful variants and of full subdiptychs, as in \textbf{Top}, when suitably adding and combining extra conditions for objects or arrows such as being Hausdorff, proper, \'{e}tale, retroconnected, and also various \emph{countability conditions}, either on the manifolds (e.g. existence of a countable dense subset) or on the maps (for instance finiteness or countability of the fibres : retro-finiteness, retro-countability).

This basic diptych is \emph{not regular}\footnote{This is indeed an important source of difficulty, but also of richness for the theory.}.

\subsubsection{Vector bundles.}
\label{vb}
Let \textbf{VecB} denote the category of (morphisms between) vector bundles (for instance in \textbf{Dif}) ; the arrows are commutative squares (see (\ref{3bt}):
\begin{diagram}[size=1.8em,tight,labelstyle=\scriptstyle,midshaft]
        E'     &\rTo^f            &         E   \\
\dSTo^{p'}         &                & \dSTo_p       \\
        B'    &\rTo^{f_0}            &         B   \\
  \end{diagram} %

There are several useful diptych structures. The basic one, denoted by {\sffamily VecB}, takes for good monos the squares with $f,f_{0}\in\cD_i$, and for good epis those with $f,f_{0}\in\cD_s$, \emph{and which moreover are $s$-full} in the sense to be defined below (\ref{3bt}). (One can also use pullback squares, which means E is \emph{induced} by $B$ along $f_0$).

\subsection{(Pre)diptych structures on simplicial (and related) categories.}
\label{simp}
Small diptychs (and especially prediptychs) may be also of interest :

\subsubsection{Finite cardinals.}
\label{ficar}
The category $\cN_c$ of all maps between \emph{finite cardinals} (or integers) defines a regular diptych $\sfNc$ by taking for $(\cN_c)_i\,/\,(\cN_c)_s$ all the injections/surjections. The set of objects is $\bbN$. The $s$-condensed objects are those which are $\neq0$.

The dual or opposite category, denoted by $\cN_c^\ast$, defines also a regular diptych $\sfNc^\ast$ with $(\cN_c^\ast)_i=((\cN_c)_s)^\ast$, $(\cN_c^\ast)_s=((\cN_c)_i)^\ast$. The product in $(\cN_c)^\ast$ is the sum in $\cN_c$. The terminal object is now $0$, and all the objects are now $s$-condensed.

\par\smallskip
When dropping $0$ one gets still diptychs denoted by $\sfNc^+$ and $\sfNc^{+\ast}$ \footnote{The latter with a terminal object, the former without such.}, the objects of which will be denoted by $\cdot n=n+1\ (n\in\bbN$)\,\footnote{The notation suggests that the elements of such an object \emph{have to be numbered from $0$ to $n$}, the dot symbolizing the added $0$.}.

\emph{The diptych $\sfNc^{+\ast}$ will be basic for our diagrammatic description of groupoids}\,\footnote{We stress again that it is \emph{not} isomorphic to $\sfNc$\,.}.

\subsubsection{Finite ordinals.}
Considering now the integers as \emph{finite ordinals}, we get the subcategories of \emph{monotone} maps which will be denoted here by $\cN_o$ , $\cN_o^+$, (denoted by $\Delta$, $\Delta^+$ in [McL]) (\emph{simplicial categories}), as well as their duals, but these have no (cartesian) products and therefore \emph{cannot define diptychs}. They define only \emph{prediptychs} denoted here by : $\sfNo$, $\sfNo^+$, $\sfNo^\ast$, $\sfNo^{+\ast}$.
\par\smallskip
\textbf{Remark.} Though $\sfNo$ is by no means isomorphic to its dual $\sfNo^\ast$, however there are two canonical isomorphisms $\Phi$, $\Psi=(\Phi^\ast)^{-1}$ (\emph{defined only on the privileged subcategories !}), which can be defined using the canonical generators as denoted in [McL]:
\begin{center}
$\Phi:(\cN_o)_i\rightarrow(\cN_o^\ast)_i=((\cN_o)_s)^\ast$, $n\mapsto\cdot n$, $\delta_j^n\mapsto(\sigma_j^{\cdot n})^\ast$,\\
$\Psi:(\cN_o)_s\rightarrow(\cN_o^\ast)_s=((\cN_o)_i)^\ast$, $\cdot n\mapsto n$, $\sigma_j^{\cdot n}\mapsto(\delta_j^n)^\ast$,
\end{center}
where a star bearing on an arrow means that this very arrow is regarded as belonging to the dual category (with source and target exchanged).

\subsubsection{Canonical prediptychs.}
\label{can}
To each (small or not) category $\cT$, we can associate \emph{three canonical prediptychs} (the last two being regular):
\begin{center}
\fbox{
$\mathsf{T}_{(\ast)}=(\cT;\cT_\ast,\cT_\ast)$\,, 
$\mathsf{T}_{(\iota)}=(\cT;\cT,\cT_{\ast})$ and $\mathsf{T}_{(\sigma)}=(\cT;\cT_{\ast},\cT)$}\,.
\end{center}
\par\smallskip

\subsubsection{Silly prediptychs.}
\label{Silpred}

Let $\un$ (or sometimes, more pictorially, $\downarrow$ or $\rightarrow$) denote the (seemingly silly) category with two objects 0, 1, and one non unit arrow 0 $\stackrel{\epsilon}{\rightarrow}$ 1 (it owns products and sums)\,\footnote{This is the category denoted by \textbf{2} in [McL], since it represents the order of the ordinal 2.}. It is canonically isomorphic with its dual $\un^{\ast}$, by exchanging the two objects. 

All the arrows are both mono- and epi-morpisms, but $\epsilon$ is not strict, and cannot be accepted as a good epi. As we shall see, it will be convenient to endow $\un$ with one of its canonical \emph{prediptych} structures (\ref{can}): 
$\un_{(\ast)}$\,,$\un_{(\iota)}$\,,$\un_{(\sigma)}$\,, according to what is needed.

\section{Commutative squares $\square$ in a diptych $\mathsf{D}=(\cD;\cD_i,\cD_s)$.}
\label{com}

\subsection{Three basic types of squares.}
\label{3bt}

Let $\square\:\cD$ \footnote{Ehresmann's notation.} denote the category of commutative squares in $\cD$ with the \emph{horizontal} composition, which can be regarded (pedantically) as the category of natural transformations between functors from $ \textsf{I}$ to $\cD$, with the (unfortunately so-called !) vertical composition. Its arrows might alternatively be described as functors from $ \textsf{I}\times \textsf{I}$ to $\cD$.

It turns out that the following properties of a commutative square $(K)$ 
\begin{diagram}[size=1.8em,tight,labelstyle=\scriptstyle,textflow]
        A'     &\rTo^{f'}            &         B'   \\
\dTo^{u}         &         {(K)}       & \dTo_{v}       \\
        A     &\rTo_{f}            &         B   \\
  \end{diagram}%
play a basic role (the terminology will be explained by the application to functors).

\par\smallskip
\textbf{Definition :} the commutative square $(K)$ is said to be :

\paragraph{}
a) \emph{$i$-faithful} if ${A'}\rITo^{(u,f')}A\times{B'}$ is a good mono ; notation :
\begin{center}
\begin{diagram}[size=2.2em,tight,labelstyle=\scriptstyle,inline]
       \ \ \ \ A'_{\SEI}     &\rTo^{f'}            &         B'   \\
\dTo^{u}         &           (K)     & \dTo_{v}       \\
        A     &\rTo_{f}            &         B   \\
\end{diagram}\ \ \ or\ \ \ 
\begin{diagram}[size=2.2em,tight,labelstyle=\scriptstyle,inline]
       A'     &\rTo^{f'}            &         B'   \\
\dTo^{u}         &           {\fid}     & \dTo_{v}       \\
        A     &\rTo_{f}            &         B   \\
\end{diagram}
\end{center}

\paragraph{}
b) \emph{a good pullback}
\footnote{In \textbf{Dif}, there are plenty of useful pullback squares generated by pairs $(f,v)$ which are not transversal in the classical sense of [B] (for instance two curves intersecting neatly in a high dimensional manifold), but there are also (actually pathological) pullback squares existing without $(f,v)$ being weakly transversal. We think the weak transversality is the most useful notion.}
 if it is a pullback square (in the usual categorical sense) and moreover $i$-faithful
\footnote{Observe that perfect squares are good pullbacks, as well as those arising from axiom (v), or more generally from (\ref{altax}), but it may happen that more general ones are needed.} ;
 one writes $f\tr v$ ($f$ and $v$ are ``\emph{weakly transversal}'') to mean that the pair $(f,v)$ can be completed in such a square ; notation :
\begin{center}
\begin{diagram}[size=2.2em,tight,labelstyle=\scriptstyle,inline]
        A'\SEpbk     &\rTo^{f'}            &         B'   \\
\dTo^{u}         &         (K)       & \dTo_{v}       \\
        A     &\rTo_{f}            &         B   \\
\end{diagram}%
\ \ \ \ or\ \ \ 
\begin{diagram}[size=2.2em,tight,labelstyle=\scriptstyle,inline]
        A'     &\rTo^{f'}            &         B'   \\
\dTo^{u}         &         {\pb}       & \dTo_{v}       \\
        A     &\rTo_{f}            &         B   \\
\end{diagram}%
\end{center}

\paragraph{}
c) \emph{$s$-full}\,\footnote{Such a square owns the \emph{parallel transfer property} : $(v\in\cD_s)\Rightarrow(u\in\cD_s)$.} if one has $f\tr v$\footnote{This was not demanded in a).} and if moreover the canonical arrow $A'\rSTo A\times_B B'$ (which is then defined) is a good epi ; notation as below.
\begin{center}
\begin{diagram}[size=2.2em,tight,labelstyle=\scriptstyle,inline]
       \ \ \ \ A'_{\SES}     &\rTo^{f'}            &         B'   \\
\dTo^{u}         &         {(K)}       & \dTo_{v}       \\
        A     &\rTo_{f}            &         B   \\
\end{diagram}\ \ \ or\ \ \ 
\begin{diagram}[size=2.2em,tight,labelstyle=\scriptstyle,inline]
       A'     &\rTo^{f'}            &         B'   \\
\dTo^{u}         &         {\ful}       & \dTo_{v}       \\
        A     &\rTo_{f}            &         B   \\
\end{diagram}
\end{center}
\subsection{Basic diptych structures on $\square\:\cD$.}
\label{DS on SQD}
The three previous kinds of squares have remarkable composition stability properties resulting from the axioms, which we shall not state here (cf Prop. A 2 of [P3], with a different terminology). The (purely diagrammatical) proofs are never very hard, but may be lengthy.

A substantial part of these properties is expressed by the following important (non exhaustive) statements, which deserve to be considered as theorems, as they bring together a very large number of various properties, which acquire much power by being gathered.

Let be given a diptych $\mathsf{D}=(\cD;\cD_i,\cD_s)$\:.

\subsubsection{Canonical structure}
\label{cast}

 On $\square\:\cD$, we get a first diptych structure , called \emph{canonical} by setting :
\begin{center}
\fbox{$\square\:\sfD=(\square\:\cD;\square\:(\cD,\cD_i),\square\:(\cD,\cD_s))$}\:,
\end{center}
which is pictured by :
\begin{center}
\begin{diagram}[size=0.9em,tight,inline,abut]%
\cdot&\rTo&\cdot\\
\dTo&&\dTo\\
\cdot&\rTo&\cdot\\
\end{diagram}
\ ;\ 
\begin{diagram}[size=0.9em,tight,inline,abut]%
\cdot&\rITo&\cdot\\
\dTo&&\dTo\\
\cdot&\rITo&\cdot\\
\end{diagram}
\ ,\ 
\begin{diagram}[size=0.9em,tight,inline,abut]%
\cdot&\rSTo&\cdot\\
\dTo&&\dTo\\
\cdot&\rSTo&\cdot\\
\end{diagram}\ \ .
\end{center}
and also, using remark (\ref{fulsub}), the full subdiptych (still defined for the diptychs of the following paragraph):
\begin{center}
\fbox{$_i\square\:\sfD=(\square\:(\cD_i,\cD);\square\:(\cD_i,\cD_i),\square\:(\cD_i,\cD_s))$}\:,
\end{center}
pictured by :
\begin{center}
\begin{diagram}[size=0.9em,tight,inline,abut]%
\cdot&\rTo&\cdot\\
\dITo&&\dITo\\
\cdot&\rTo&\cdot\\
\end{diagram}
\ ;\ 
\begin{diagram}[size=0.9em,tight,inline,abut]%
\cdot&\rITo&\cdot\\
\dITo&&\dITo\\
\cdot&\rITo&\cdot\\
\end{diagram}
\ ,\ 
\begin{diagram}[size=0.9em,tight,inline,abut]%
\cdot&\rSTo&\cdot\\
\dITo&&\dITo\\
\cdot&\rSTo&\cdot\\
\end{diagram}\ \ .
\end{center}
Defining an analogous full subdiptych with $\cD_s$ replacing $\cD_i$ \emph{requires a more restrictive choice for the squares taken as good epis}. (These full subdiptychs will be essential to get fibred products of diagrams, hence of groupoids).

\subsubsection{Full and pullback epis.}
\label{fpe}
We have the two basic diptychs:

\begin{center}
\fbox{$\square_{\mathsf{(i,ful)}}\:\sfD=%
(\square\:\cD;\square\:(\cD,\cD_i),\ful(\cD,\cD_s))$}\\
\fbox{$\square_{\mathsf{(i,pb)}}\:\sfD=%
(\square\:\cD;\square\:(\cD,\cD_i),\pb(\cD,\cD_s))$}
\end{center}
which can be pictured by :
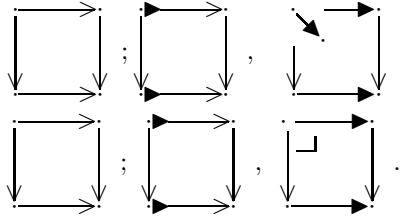
\begin{center}
\begin{diagram}[size=1.6em,tight,inline,abut]%
\cdot&\rTo&\cdot\\
\dTo&&\dTo\\
\cdot&\rTo&\cdot\\
\end{diagram}
\ ;\ 
\begin{diagram}[size=1.6em,tight,inline,abut]%
\cdot&\rITo&\cdot\\
\dTo&&\dTo\\
\cdot&\rITo&\cdot\\
\end{diagram}
\ ,\ 
\begin{diagram}[size=1.6em,tight,inline,abut]%
\ \ \ \cdot_{\SES}&\rSTo&\cdot\\
\dTo&&\dTo\\
\cdot&\rSTo&\cdot\\
\end{diagram}\
\end{center}

\begin{center}\ \ 
\begin{diagram}[size=1.6em,tight,inline,abut]%
\cdot&\rTo&\cdot\\
\dTo&&\dTo\\
\cdot&\rTo&\cdot\\
\end{diagram}
\ ;\ \ 
\begin{diagram}[size=1.6em,tight,inline,abut]%
\cdot&\rITo&\cdot\\
\dTo&&\dTo\\
\cdot&\rITo&\cdot\\
\end{diagram}
\ ,\ \ 
\begin{diagram}[size=1.6em,tight,inline,abut]%
\cdot_{\sepb}&\rSTo&\cdot\\
\dTo&&\dTo\\
\cdot&\rSTo&\cdot\\
\end{diagram}\ \ .
\end{center}
and in this way, thanks to a parallel transfer property, we get the expected \emph{basic full subdiptychs} :

\begin{center}
\fbox{$_s\square_{\mathsf{(i,ful)}}\:\sfD=%
(\square\:(\cD_s,\cD);\square\:(\cD_s,\cD_i),\ful(\cD_s,\cD_s))$}\\
\fbox{$_s\square_{\mathsf{(i,pb)}}\:\sfD=%
(\square\:(\cD_s,\cD);\square\:(\cD_s,\cD_i),\pb(\cD_s,\cD_s))$}
\end{center}
pictured by :
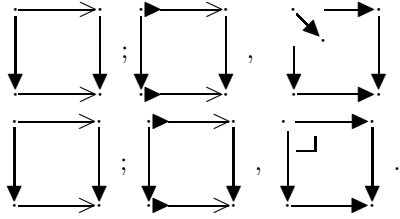
\begin{center}
\begin{diagram}[size=1.6em,tight,inline,abut]%
\cdot&\rTo&\cdot\\
\dSTo&&\dSTo\\
\cdot&\rTo&\cdot\\
\end{diagram}
\ ;\ 
\begin{diagram}[size=1.6em,tight,inline,abut]%
\cdot&\rITo&\cdot\\
\dSTo&&\dSTo\\
\cdot&\rITo&\cdot\\
\end{diagram}
\ ,\ 
\begin{diagram}[size=1.6em,tight,inline,abut]%
\ \ \ \cdot_{\SES}&\rSTo&\cdot\\
\dSTo&&\dSTo\\
\cdot&\rSTo&\cdot\\
\end{diagram}\
\end{center}

\begin{center}\ \ 
\begin{diagram}[size=1.6em,tight,inline,abut]%
\cdot&\rTo&\cdot\\
\dSTo&&\dSTo\\
\cdot&\rTo&\cdot\\
\end{diagram}
\ ;\ \ 
\begin{diagram}[size=1.6em,tight,inline,abut]%
\cdot&\rITo&\cdot\\
\dSTo&&\dSTo\\
\cdot&\rITo&\cdot\\
\end{diagram}
\ ,\ \ 
\begin{diagram}[size=1.6em,tight,inline,abut]%
\cdot_{\sepb}&\rSTo&\cdot\\
\dSTo&&\dSTo\\
\cdot&\rSTo&\cdot\\
\end{diagram}\ \ .
\end{center}

\par\bigskip
One also defines two basic \emph{subdiptychs} by taking :
\begin{center}
\fbox{$\pb\:\sfD=(\pb\:\cD;\pb\:(\cD,\cD_i),\pb\:(\cD,\cD_s))$}
\end{center}

\begin{center}
\fbox{$_s\ful\:\sfD=(\ful\:(\cD_s,\cD);\ful\:(\cD_s,\cD_i),\ful\:(\cD_s,\cD_s))$}
\end{center}

pictured by :

\begin{center}\ \ 
\begin{diagram}[size=1.2em,tight,inline,abut]%
\cdot_{\sepb}&\rTo&\cdot\\
\dTo&&\dTo\\
\cdot&\rTo&\cdot\\
\end{diagram}
\ \ ;\ \ 
\begin{diagram}[size=1.2em,tight,inline,abut]%
\cdot_{\sepb}&\rITo&\cdot\\
\dTo&&\dTo\\
\cdot&\rITo&\cdot\\
\end{diagram}
\ ,\ \ \ 
\begin{diagram}[size=1.2em,tight,inline,abut]%
\cdot_{\sepb}&\rSTo&\cdot\\
\dTo&&\dTo\\
\cdot&\rSTo&\cdot\\
\end{diagram}
\end{center}

\begin{center}\ \ 
\begin{diagram}[size=1.2em,tight,inline,abut]%
\cdot&\rTo&\cdot\\
\dSTo&\ful&\dSTo\\
\cdot&\rTo&\cdot\\
\end{diagram}
\ \ ;\ \ 
\begin{diagram}[size=1.2em,tight,inline,abut]%
\cdot&\rITo&\cdot\\
\dSTo&\ful&\dSTo\\
\cdot&\rITo&\cdot\\
\end{diagram}
\ ,\ \ \ 
\begin{diagram}[size=1.2em,tight,inline,abut]%
\cdot&\rSTo&\cdot\\
\dSTo&\ful&\dSTo\\
\cdot&\rSTo&\cdot\\
\end{diagram}.
\end{center}

\subsubsection{Iteration.}
\label{iter}
One immediately notices than this theorem allows an \emph{iteration} of the construction of commutative squares giving rise to new diptych structures for commutative cubes, and so on, which would be very difficult to handle directly.
\par\smallskip
This is especially interesting when dealing with commutative cubes since such a diagram gives rise to \emph{three commutative squares} in $\square\:\cD$, each edge of which (which is actually a square of $\cD$) belonging to \emph{two different squares} of $\square\:\cD$, and this gives a powerful method for deducing properties of certain faces (for instance being a pullback), or edges from properties of the others, using diptych properties of parallel transfer. But we cannot develop more here.

\section{Diagrams of type $\sfT$ in a diptych.}
\label{typeT}
We are now going to replace the silly category $ \textsf{I}$ of \ref{Silpred} by a more general notion generalizing the construction of commutative squares, and allowing to perform various set-theoretic constructions in a general diptych $\sfD$.

\par\medskip
Let $\sf{T}$ be a \emph{small} prediptych\,\footnote{It may be sometimes useful to extend the following definitions when $\sf{T}$ is just a graph with two given subgraphs.} (\ref{predip}), and let $\mathsf{D}=(\cD;\cD_i,\cD_s)$ be a diptych.

\subsection{Definitions : objects and morphisms.}
\label{defdiag}

We denote by \fbox{$\square^{\cT}\:\cD:\,=\cD^{\cT}$} the category having :
\par\smallskip
-- as \emph{objects} the elements of 
\fbox{$\un\,^{\sfT}\:\sfD:\,=\Hom(\sfT,\sfD)$}
\,, called \emph{diagrams of $\sfD$ of type $\sfT$}, which means those \emph{functors} $F$ from $\cT$ to $\cD$ such that one has :
$$F(\cT_i)\subset\cD_i \text{and} F(\cT_s)\subset\cD_s\leqno\quad(1)$$
(these last conditions are void when $\sfT=\sfT_{(\ast)}$(\ref{can}))\,\footnote{Forgetting the prediptych structure of $\sfT$ and dropping conditions (1) would oversimplify the theory, but deprive it of all its strength.} ;
\par\smallskip
-- as \emph{arrows} the \emph{morphisms} (i.e. the natural transformations) between such functors\footnote{Again with the ``vertical composition'', which we prefer here to write horizontally, \emph{drawing the diagrams vertically}, and the morphisms horizontally. Of course one can exchange everywhere \emph{simultaneously} ``vertical'' and ``horizontal'', since the distinction is purely psychological and notational, and since $\square^{\mathsf{vert}}\:\cD$ and $\square^{\mathsf{hor}}\:\cD$ are canonically isomorphic.}.

 These morphisms may be described :
\par\smallskip
- either as functors $\Phi$ from $\cT\times \un$ to $\cD$ such that (denoting by $(:)$ the set of the two units of $ \un$):
\begin{center}
$\Phi(\cT_i\times(:))\subset\cD_i$ and $\Phi(\cT_s\times(:))\subset\cD_s$
\end{center}
\par\smallskip
- or as functors $\Psi$ from $\cT$ to $\square^{\mathsf{vert}}\:\cD$\,
\footnote{This notation means that $\square\:\cD$ has to be considered here with its vertical composition law.}
such that, with the notations of (\ref{DS on SQD}) :
\begin{center}
$\Psi(\cT_i)\subset\square\:(\cD_i,\cD)=\ _i\square\:\cD$ and $\Psi(\cT_s)\subset\square\:(\cD_s,\cD)=\ _s\square\:\cD$\ .
\end{center}
In words this means that such a morphism may be viewed :
\par\smallskip
- \emph{either as a diagram in} $\sfD$\emph{ of type} %
$\sfT\times\, \un_{(\ast)}$ (cf.\ref{can} and \ref{Silpred})
\par\smallskip
- \emph{or as a diagram of the same type} $\sfT$ \emph{in} $\square^{\mathsf{vert}}\:\sfD$, regarded with its \emph{canonical} (vertical) diptych structure (\ref{DS on SQD}).

\par\smallskip
 However be careful that, in such a description, though the previous \emph{definition} of the morphisms uses the \emph{vertical} composition of squares, the \emph{composition} of these morphisms (called vertical in [McL]!!) involves the \emph{horizontal} composition of squares.

\subsection{Diptych structures on the category of diagrams.}
\label{dipdiag}

\subsubsection{Definitions.}
\label{defdipdiag}

Using the latter interpretation, and taking now into account the \emph{horizontal} composition, it is clear that any prediptych structure on $\square\:\cD$ determines a prediptych structure on $\square^{\cT}\:\cD$. For instance, using on $\sfT$ its trivial prediptych structure, we can consider the canonical \emph{prediptych} structure $\square^{\sfT}\:\sfD$, but this \emph{not} in general a dyptich structure (one needs parallel transfer properties for good epis in order to ensure conditions (1) of (\ref{defdiag})).

Using the parallel transfer properties of $s$-full and pullback squares, one can get three useful \emph{diptych structures} denoted by:

\begin{center}
\fbox{
$(\square_{\mathsf{(i,ful)}})^{\sfT}\:\sfD$\,, 
$(\square_{\mathsf{(i,pb)}})^{\sfT}\:\sfD$ \, and $\pb^\sfT\:\sfD$}\,.
\end{center}

Of special interest, as we shall see, will be those diagrams which preserve certain pullbacks, since groupoids may be described by diagrams of this type.

\subsubsection{The silly case.}
\label{silcase}

\emph{The case of commutative squares is just the special case when one takes for} 
$\sfT$ \emph{the silly category} $\un$\,, \emph{with one of its prediptych structures} (\ref{Silpred}), since one has:
\begin{center}
$\un\,^{\un_{(\ast)}}\:\cD=\left|\cD\right|$ , 
$\un\,^{\un_{(\iota)}}\:\cD=\left|\cD_i\right|$ , 
$\un\,^{\un_{(\sigma)}}\:\cD=\left|\cD_s\right|$ , \footnote{In the second members of this line, the signs $\left|\ \right|$ denote the forgetful functor forgetting the composition laws, since they are not defined by the first members.}\\
$\square\,^{\un_{(\ast)}}\:\cD=\square\:\cD$ , 
$\square\,^{\un_{(\iota)}}\:\cD={_i\square}\:\cD$ ; 
$\square\,^{\un_{(\sigma)}}\:\cD={_s\square}\:\cD$ ; 
\end{center}
more precisely :

\begin{center}
\fbox{%
$\square\,^{\un_{(\ast)}}\:\sfD=\square\:\sfD$ , 
$\square\,^{\un_{(\iota)}}\:\sfD=\,_i\square\:\sfD$
}\\
\fbox{
$(\square_{\mathsf{(i,ful)}})^{\un_{(\ast)}}\:\sfD=
\square_{\mathsf{(i,ful)}}\:\sfD$ , 
$(\square_{\mathsf{(i,pb)}})^{\un_{(\ast)}}\:\sfD=
\square_{\mathsf{(i,pb)}}\:\sfD$\:.
}%
\end{center}
and also :
\begin{center}
\fbox{
$(_s\square_{\mathsf{(i,ful)}})^{\un_{(\sigma)}}\:\sfD=\:
_s\square_{\mathsf{(i,ful)}}\:\sfD$ , 
$(_s\square_{\mathsf{(i,pb)}})^{\un_{(\sigma)}}\:\sfD=\:
_s\square_{\mathsf{(i,pb)}}\:\sfD$\:.
}%
\end{center}

\subsection{The exponential law for diagrams.}
\label{exp}
Given another similar datum $\mathsf S$, we have the \emph{exponential law}, which, with our notations, reads (at least when $\sfS$ and $\sfT$ are trivial diptychs) :
\begin{center}
\fbox{
$\square^{\sfS}(\square^{\sfT}\:{\sfD})=\square^{\sfS\times\sfT}\:\sfD$}\, 
\end{center}
(with a lot of possible variants), and the canonical isomorphism : $\sfS\times\sfT\approx\sfT\times\sfS$ yields a canonical isomorphism :
\begin{center}
\fbox{$\square^{\sfS}(\square^{\sfT}\:{\sfD})\approx\square^{\sfT}(\square^{\sfS}\:{\sfD})$}\,.
\end{center}

Particularly one has: 
\fbox{$(\square\,\sfT)^n=\square^{\sfT^n}$}, \\
and, specializing for $\sfT=\un_{(\ast)}$: 
\fbox{$(\square)^n=\square^{(\un_{(\ast)})^n}$}\,.
\par\bigskip
\textbf{Remarks about notations.} The reader may have observed that the symbols $\square$ and $\un$ used above are intendedly ambiguous and protean, since this allows to memorize and visualize a lot of various properties :

-- the symbols $\un$ or $\un_{(\bullet)}$, where $\bullet$ stands for $\ast$, $\iota$, or $\sigma$ :\\
- denote the silly category, possibly with various prediptych structures on it;\\
- denote the bifunctor Hom\,(?,?), with the first argument treated as an exponent);\\
- often behaves formally as a 1;\\
- may sometimes suggest the functor $\sfT\mapsto\sfT\times\un_{(\bullet)}$ ;

-- the symbol $\square$ with possibly labels in various positions $_\diamond\square\text{\tiny{\hskip-5.5pt{?}}}\phantom{a}_{(\sharp,\natural)}^{\top}$ :\\
- may picture or suggest the product $\un\times\un$ or more precisely various instances of $\un_{(\bigtriangleup)}\times\un_{(\bigtriangledown)}$\,;\\
- may create the commutative squares of a category or of a diptych with possible extra conditions :

$\ast$ on the left, they bear on the vertical edges ;

$\ast$ on the right, they describe the subcategories for a (pre)diptych structure ;

$\ast$ inside, they describe global properties of the squares ;\\
- may become the bifunctor $\mathcal H$om\,(?,?), with possible restrictions on the squares involved.

\section{Groupoids as diagrams in $\sfE=\mathsf{Set}$.}
\label{strgrpd}
We shall be concerned only with \emph{small} groupoids, viewed as generalizing both groups and graphs of equivalence relations, as well as specializing (small) categories.

According to our program, we need a diagrammatic description of both groupoid data and axioms, which will be transferred from the diptych $\sfE=\mathsf{Set}$ to a general diptych $\sfD$.

This can be achieved in two complementary ways :

- either using finite diagrams (called \emph{sketches} in Ehresmann's terminology) ;

- or using the \emph{simplicial} description by the \emph{nerve}.

Though seemingly more abstract, the latter turns out to be often the most convenient for theoretical purposes, while the former is adapted to practical handling.

We refer to [McL] and also to (\ref{simp}) for general properties and notations concerning simplicial categories and \emph{simplicial objects}, which are just functors from $\cN_o^{+\ast}$ to $\cE$, and more precisely \emph{diagrams of type $\sfNo^{+\ast}$ in $\sfE$}. 

\subsection{Some remarks about $\sfNc^+$ and $\sfNo^+$.}(see \ref{simp})
\subsubsection{The arrows as functors.}
\label{arf}
For each $n\in\bbN$, we shall denote by :\\
\vspace{0.5pt}
$\ast$ $\overrightarrow{\cdot n}$ the ordinal $\cdot n=n+1$ viewed as a small \emph{category} (defining the order)\,\footnote{$\,\overrightarrow{\cdot 1}\,=\,\overrightarrow{2}$ is just what we called in (\ref{Silpred}) the silly category $\un=\,\downarrow\,=\,\rightarrow\,$.};\\
\vspace{0.5pt}
$\ast$ $\dot{\wedge} n$ the subcategory consisting of arrows of this order category with source $0$\footnote{And of course the units.};\\
\vspace{0.5pt}
$\ast$ $\overleftrightarrow{\cdot n}=(\cdot n)^2$ the banal \emph{groupoid} $(\cdot n)\times(\cdot n)$\,, which ignores the ordering\,\footnote{Though the numbering of the base is kept!}.

With these structures on the objects, the arrows of $\cN_o^+$ and $\cN_c^+$ may be regarded (though this looks pedantic) as \emph{functors (or morphisms) between small categories}.

It may be sometimes suggestive to write :
\begin{center}
$\overleftrightarrow{\cdot 0}=\,\overleftrightarrow{1}=\cdot$\,, $\overleftrightarrow{\cdot 1}=\,\overleftrightarrow{2}=\,\leftrightarrow$\,, 
$\dot{\wedge} 2=\wedge$\,,
$\overleftrightarrow{\cdot 2}=\,\overleftrightarrow{3}=\bigtriangleup$\,,
$\overleftrightarrow{\cdot 3}=\,\overleftrightarrow{4}=\boxtimes$\,.
\end{center}

\subsubsection{Pullbacks in $\cN_c^\ast$.}\label{pbn}
Pullbacks in $\cN_c^\ast$ come from pushouts in $\cN_c$\,.
One can check that \emph{all the commutative squares describing the relations between the canonical generators of} $\cN_o$ (see [McL]) \emph{become pushouts} 
\emph{when written in} $\cN_c$, \emph{and they generate, by composition, all the pushouts of} $\cN_c$. More precisely (with the notations in [McL]) the squares involving the injections $\delta_j$\,'s alone or mixing both $\delta_j$\,'s and $\sigma_k$\,'s \emph{are pullbacks too}, but \emph{not those with  the surjections $\sigma_j$\,'s alone}\,\footnote{This induces a strong dissymmetry between the dual diptychs $\sfNc$ and $\sfNc^\ast$.}.

This is still valid when dropping $0$ (but one looses the squares describing the sums as pushouts).

\subsection{Characterization of the nerve of a groupoid.}
\label{charnerv}
The special properties of a groupoid among categories yield very special and remarkable properties as well as alternative descriptions for its nerve.

\subsubsection{Three descriptions for the nerve.}
\label{desner}
Given a (small) groupoid $\mathbf G$, denoted loosely by $G$ or $G\rightrightarrows B$, we can associate to it \emph{three canonically isomorphic simplicial objects} (among which the first one is the nerve of $G$ regarded as a category).
First we define and denote the three images of the generic object $\cdot n=n+1$\,\footnote{Though these images may be regarded as instances of diagrams of finite type in $G$, the notations used below differ slightly from those used in the previous section for the general case.}:
\begin{enumerate}
	\item $^{\downarrow}G^{(n)}=\hom(\overrightarrow{\cdot n},G)=\left\{\text{paths of }G\text{ of length }n\right\}$ ;
	\vspace{0.5pt}
	\item $^{\wedge}G^{(n)}=\hom(\dot{\wedge} n,G)=\left\{n\,\text{-uples of arrows of }G\text{ with the same source}\right\}$;
	\vspace{0.5pt}
	\item $^{\updownarrow}G^{(n)}=\hom(\overleftrightarrow{\cdot n},G)=\left\{\text{commutative diagrams of }G\,\text{with }{\cdot n}\,\text{vertices}\right\}$\footnote{With two-sided edges and numbered vertices.}.
\end{enumerate}
\emph{Identifying} these three sets, we can write loosely:
$$G^{(n)}=\,^{\downarrow}G^{(n)}=\,^{\wedge}G^{(n)}=\,^{\updownarrow}G^{(n)}\,.$$
We shall also sometimes feel free to write loosely and suggestively :
\begin{center}
$G^{(0)}=B$\,, $G^{(1)}=\,^{\downarrow}G=\,^{\updownarrow}G=G$\,, $G^{(2)}=\,\wedge G=\bigtriangleup G$\,, $G^{(3)}=\boxtimes\,G$\,.
\end{center}
The images of the canonical generators of $\cN_c^+$ are then especially easy to define with the interpretation (3), since they consist in repeating or forgetting one vertex.

But, with the descriptions (1) and (3), one can immediately interpret the required contravariant functors from $\cN_o^+$ to $\cE$ as being just special instances of the classical contravariant hom-functors $\hom(?,G)$, which consist in letting the arrows of $\cN_o^+$ act \emph{by right morphism composition} (in the category of morphisms or functors between small categories) with the elements of $^{\downarrow}G^{(n)}$ or $^{\updownarrow}G^{(n)}$.

\subsubsection{Properties of the nerve of a groupoid.}\label{proner}
\paragraph{}
- (\emph{Extension of the nerve}).
With the interpretation (3), we get a bonus, since it is now obvious that this contravariant functor \emph{extends to $\cN_c^+$}, and so defines a diagram in $G$ of type $\sfNc^{+\ast}$ . (We remind that, in this interpretation, integers are interpreted as small banal groupoids, and the arrows of $\cN_c^+$ as morphisms.)
\paragraph{}
- (\emph{Exactness properties}).
Moreover this extended functor, which is now defined on a \emph{diptych}, is ``\emph{exact}'', there meaning that \emph{it preserves pullbacks}\,\footnote{But \emph{not products} (and not the pushouts which don't arise from perfect squares). When some risk of confusion might arise,
 it would be more correct to say something like ``diptych-exact'' or ``p.b.-exact'', since here the term ``exact'' has to be understood in a much weaker sense than the general meaning for functors between categories : it has not to preserve all (co)limits.}. It would be enough to check this for the generating pullbacks of $\cN_c^\ast$ (\ref{pbn}).
\paragraph{}
- \emph{Conversely} one might check that the datum of an exact diptych morphism $\mathbf G=(G^{(n)})_{(n\in\bbN)}$\,\footnote{Relaxed notation, using only the images of the objects.} from $\sfNc^{+\ast}$ to $\sfE=\mathsf{Set}$ determines a groupoid, the nerve of which is the restriction of this morphism to $\sfNo^{+\ast}$.
\paragraph{}
- (\emph{Sketch of a groupoid}).\,\footnote{This notion was introduced by Ehresmann for general categories.} Actually one might check that the groupoid \emph{data} are fully determined by the restriction of $\mathbf G$ to the (truncated) full subcategory $[\cdot\textbf{\large 2}]$ of $\cN_c^{+\ast}$ generated by $\cdot0,\cdot1,\cdot2$, while the groupoid \emph{axioms} are expressed by its restriction to $[\cdot\textbf{\large3}]$ (and the conditions that the images of the previous generating squares be pullbacks)\,\footnote{We cannot give more details here.}.

\subsubsection{Concrete description.}
\label{codes}
Actually the previous data are somewhat redundant, and it is convenient for our purpose to observe that $\mathbf G$ may be fully described by the data $(G,B,\omega_G,\alpha_G,\delta_G)$\,(satisfying axioms which we shall not make explicit), where :
\begin{itemize}
	\item $B$ (``base'') and $G$ (``set of arrows'') are objects of $\cE$ ;
	\item $\omega_G:B\rITo G$ is an injection (``unit law'') ;
	\item $\alpha_G:G\rSTo B$ is a surjection (source map) ;
	\item $\delta_G:\wedge G=\bigtriangleup G\rSTo G$ is a surjection (``division map''\footnote{This means the map : $(y,x)\mapsto yx^{-1}$.}), 
with :\\
$\wedge G=G\times_B G$ (fibred product of $\alpha$ by $\alpha$).
\end{itemize}

One often writes $\alpha$ for $\omega\alpha$, and we set : 
\begin{itemize}
	\item 
$\tau_G=(\beta_G,\alpha_G):G\rightarrow B\times B$\quad
 (``anchor map''\,\footnote{Mackenzie's terminology.}, or ``\emph{transitor}'').
\end{itemize}

From these data, it is not difficult to get the range map $\beta_G$, the inverse law $\iota_G$, and the composition law $\gamma_G$ (defined on the fibred product of $\alpha_G$ and $\beta_G$).

With notations from Mac Lane:
$\omega_G$, $\alpha_G$, $\beta_G$, $\delta_G$
 would be the respective images of :
$(\sigma_0^1)^\ast,\,(\delta_1^1)^\ast,\,(\delta_0^1)^\ast,(\delta_0^2)^\ast\,$ ( in the dual $\sfNc^\ast$).

\section{$\sfD$-groupoids.}

\subsection{Definition of $\sfD$-groupoids.}

\subsubsection{Groupoids in a diptych.}
\label{groudi}

The diagrammatic description of a set-theoretical groupoid given above leads to \emph{define} a $\sfD$-groupoid as an \emph{exact diptych morhism} (\ref{proner}) , i.e. \emph{preserving good monos, good epis, and good pullbacks}:
$$\mathbf {G}:\sfNc^{+\ast}\rightarrow\sfD\,.$$

As above, the exactness properties allow to characterize $\mathbf G$ by restriction to various subcategories, and to recover in this way the simplicial description as well as the skecth ones and the $(G,B,\omega_G,\alpha_G,\delta_G)$ presentation.

We shall also use as previously the relaxed notations $(G^{(n)})_{n\in\bbN}$ (omitting the effect of $\mathbf G$ on the arrows), or briefly $G$ (meaning $G^{(1)}$), $B$ for $G^{(0)}$, and $\gPdS{G}{B}$ and its variants.

We shall denote by $\rm{Gpd}(\sfD)$ the class of $\sfD$-groupoids.
\subsubsection{Examples.}
This very general notion may be specialized using for instance the various examples given in (\ref{dipex}).

Note that the \emph{$\mathsf{Top}$\,-groupoids} are those for which the source map has to be \emph{open}, a condition which can hardly be avoided for getting a useful theory.

\emph{Lie groupoids}\,\footnote{Introduced by C. Ehresmann under the name of \emph{differentiable groupoids}.} are of course the $\mathsf{Dif}$\,-groupoids, with the usual diptych structure on $\mathbf{Dif}$, but, using some of the above-mentioned variants, the theory will include, for instance, among others, étale or $\alpha$\,-connected groupoids as well.

{\sffamily VecB}-\emph{groupoids} were used in [P6].

\subsubsection{Null groupoids.}
\label{nulg}
 For any object $B$ of $\cD$, the constant simplicial object $B$ is a $\sfD$-groupoid, called ``\emph{null}''\,\footnote{We cannot accept the traditional terminology ``discrete'', which has another topological meaning, and which moreover does not agree with the group terminology, unlike ours.}, and denoted by $\nul B$. All arrows are units, and the map $\omega_G$ is an isomorphism.

\subsection{Principal and Godement $\sfD$-groupoids.}
\label{ppal}

\subsubsection{The transitor $\tau_G$.}
From a purely set-theoretical point of view, the transitor 
$\tau_G:G\rightarrow B\times B$
measures the \emph{(in)transitivity}\,\footnote{Whence the terminology adopted here, preferably to ``anchor map''.} of the groupoid $G$.

It turns out that, in the diptych setting, its properties encapsulate a very rich ``structured'' information. We just mention, without developing here, that, for instance in \textbf{Dif}, one can fully characterize, just by very simple properties of $\tau_G$, not only the graphs of regular equivalences, but (among others) the gauge groupoids of principal bundles, the Poincaré groupoids of Galois coverings, the holonomy groupoids of foliations, the Barre Q-manifolds, the Satake V-manifolds (or orbifolds).

\subsubsection{Principal $\sfD$-groupoids.}
\label{palg}
The notion of \emph{graph of regular equivalence} (\ref{Regeq}) may be carried over in any diptych as follows.

Given a good epi $B\rSTo^{q\quad} Q$ of $\sfD$, we can construct the iterated fibred product of $q$, denoted by 
$R^{(n)}=\,{\scriptstyle\stackrel{\cdot n}{\times}_Q}B$ $(n\in\bbN)$
and check that 
$(n\mapsto R^{(n)})_{(n\in\bbN)}$ allows to define a $\sfD$-simplicial object which is a groupoid. We shall say that $\text{\textbf R}=(R^{(n)})_{n\in\bbN}$ is the \emph{principal groupoid associated to $q$} (with base $B$).

Moreover we have an ``\emph{augmentation}'', which means an extended ``exact'' (\ref{proner}) diptych morphism (\ref{proner}) from $\sfNc^\ast$ to $\sfD$. This augmentation carries $0$\,\footnote{Which might be written $0={\cdot(-1)\,}$.} to $R^{(-1)}=Q$, and the added generator $\eta^{\ast}=(\delta_0^0)^{\ast}$\,\footnote{Dual arrow, from Mac Lane's notations.} to $B\rSTo^q Q$.

This situation gives rise to a perfect square (\ref{altax}):
\begin{diagram}[size=2em,tight,labelstyle=\scriptstyle,textflow]
        R     &\rSTo^{\beta_R}            &         B   \\
\dSTo^{\alpha_R}         &            \pile{\pb \\ \po}    & \dSTo_{q}       \\
        B     &\rSTo_{q}            &         Q   \\
  \end{diagram}

The pushout property of this square shows that the good epi $q$ is \emph{uniquely determined} by the knowledge of the $\sfD$-groupoid $R$.

Applied to $G\rSTo^{\alpha_G} B$, this construction allows considering $\wedge G$ as a principal groupoid (see below (\ref{canac})).

\subsubsection{Banal groupoids.}
\label{bagr}
This construction applies in particular when $B$ is an \emph{$s$-condensed object} (\ref{TO}). Then one has $R^{(n)}=B^{\cdot n}$ ($\cdot n$ times iterated product). It is called the ``\emph{banal groupoid}''\,\footnote{We cannot accept the term ``coarse'', often used in the literature for the same reason as for ``discrete''.} associated to $B$ ; it is principal. It is denoted by $B\times B$ or $B^2$.

More generally such a banal groupoid (possibly non principal in the absence of a terminal object) is associated to an object $B$ whenever the canonical projections ${\rm{pr}}_i:B\times B\rightarrow B$ are good epis\,\footnote{We gave above examples of diptychs in which such is not always the case.}.

One might call ``\emph{proper}'' those objects $B$ for which the banal groupoid $B^2$ is defined. When such is the case, the transitor $\tau_G$ may be viewed as a $\sfD$-functor.
\par\smallskip

Given an integer $n$ , we can attach to it, besides (for $n\neq0$) the banal groupoid $\overleftrightarrow n =n^2$ (which uses the product in $\cN_c$), a banal cogroupoid $^2n$ or $2n$, using the product $\times^\ast$ of $\cN_c^\ast$ (which is the coproduct $+$ of $\cN_c$).

\subsubsection{Godement groupoids and Godement diptychs.}
\label{God}%

\textbf{Definition.} A $\sfD$-groupoid $G$ is called a \emph{Godement groupoid} if $\tau_G:G\rITo B\times B$ is a good mono.

Every principal groupoid is a Godement groupoid.

$\bullet$ We say the \emph{Godement axiom} is fulfilled, and $\sfD$ is a \emph{Godement diptych} if \emph{conversely} every Godement $\sfD$-groupoid is principal.
\par\smallskip
The fact that the diptych {\sffamily Dif} (\ref{man}) is Godement is the content of the so-called \emph{Godement theorem}, proved in Serre [LALG]\,\footnote{Where the formal aspect of this theorem is clearly visible, and inspired our Godement axiom.}.

But it is highly remarkable that nearly all of the examples of diptychs given above are indeed Godement diptychs, as well as most of the diagram diptychs we constructed above, provided one starts with a Godement diptych. Such a statement includes a long list of theorems, which are not always classical.

$\bullet$ \emph{From now on, we shall assume $\sfD$ is Godement whenever this is useful.}

\subsection{Regular groupoids.}
More generally, a $\sfD$-groupoid $G$ is called \emph{regular} if $\tau_G$ is regular (\ref{Dipda}). Then we have the following factorization of $\tau_G$:
$$G\rSTo^\pi R\rITo^{\tau_R} B\times B\,,$$
where R is a Godement $\sfD$-groupoid, hence principal, so that we can construct the perfect square (\ref{palg}) ; then the orbit space $Q$ exists as an object of $\sfD$, and is also the pushout of $\alpha_G$ and $\beta_G$.

However the object $Q$ inherits an ``\emph{extra structure}'' from the arrow $G\rSTo^\pi R$, a good epi which \emph{measures how much the $s$-full pushout square $GBBQ$, fails to be a pullback}.

\par\smallskip
$\bullet$ \emph{In the general case, when $G$ is neither principal nor even regular, the aim of the present paper is to define a kind of ``virtual augmentation'' (as a substitute for the failing one), which is the $\sfD$-Morita equivalence class of $G$, and which has to be considered intuitively as defining the ``\emph{virtual structure}'' of the orbit space.}

\subsubsection{Plurigroups.}
An important special case of regular $\sfD$-groupoid is when $R=\nul B$ (\ref{nulg}) (this is indeed equivalent to $\alpha_G=\beta_G$). When such is the case, we shall say $G$ is a ``$\sfD$-\emph{plurigroup}''\,\footnote{We keep the term \emph{$\sfD$-group} for the case when moreover $B$ is a \emph{terminal object}. On the other hand, the possible term ``multigroup'' would create confusion with the multiple categories, which have nothing to do with the present case.}.

\subsubsection{$s$-transitive $\sfD$-groupoids.}
The opposite degeneracy is when $\tau_G$ is a good epi: we shall say $G$
is $s$-transitive. 

When $\sfD=\text{\sffamily Set}$, this just means that the orbit space is reduced to a singleton, but in {\sffamily Dif}, for instance, this has very strong implications, since this means essentially that $G$ may be viewed as the gauge groupoid of a principal bundle\,\footnote{The term ``Lie groupoid'' was first reserved to that special case (see for instance the first textbook by K. Mackenzie), till A. Weinstein and P. Dazord changed the terminology, with my full agreement.}. This can be proved in a purely diagrammatic way, which allow to extend these concepts to all (Godement) diptychs. In fact the orbit space has to be thought as ``a singleton structured by a group''.

\par\medskip
One of the basic reason of the strength of the notion of $\sfD$-groupoid is that it unifies and gathers in a single theory all these various degeneracies.

\section{The category $\mathbf{Gpd}(\sfD)$.}

\subsection{$\sfD$-functors.}
\subsubsection{$\sfD$-functors as natural transformations.}
\label{Dfun}
\emph{$\sfD$-functors} (or \emph{morphisms}) between $\sfD$\,-groupoids are of course special cases of diagram morphisms (\ref{typeT}) and, as such, are defined as natural transformations between the $\sfD$\,-groupoids viewed as functors.

But the pullback property of the generating squares of the diagrams defining groupoids (\ref{proner}) has very strong implications (arising from the preliminary study of commutative squares in a diptych) which we cannot develop here, referring to [P4] for more details. We just mention a few basic facts.
\par\smallskip
A $\sfD$\,-functor $\mathbf{f : H\rightarrow \mathbf G}$ is fully determined by $f^{(1)}:H^{(1)}\rightarrow G^{(1)}$ and hence often denoted loosely by $f:H\rightarrow G$ ; we shall write : $H^{(0)}=E$, $ G^{(0)}=B$. It is called \emph{principal} if $H $ is principal.

We say $f$ is a $i$/$s$-functor when $f$ lies in $\cD_i$/$\cD_s$ ; as a consequence, one can check this is still valid for all the $f^{(n)}$.

We get in this way the category $\mathbf{Gpd}(\sfD)$ , with $\rm{Gpd}(\sfD)$ as its base, of ($\sfD$-functors or morphisms between) $\sfD$-groupoids, and two subcategories $\mathbf{Gpd}_i(\sfD)$ and $\mathbf{Gpd}_s(\sfD)$, but, as announced earlier, the second one will \emph{not} be the right candidate for good epis in $\mathbf{Gpd}(\sfD)$ (see below).

An arrow of $\mathbf{Gpd}(\sfD)$ is said to be \emph{split} if it is right invertible, in other words if it admits of a section.

\subsubsection{$\sfD$-functors as $(\square^{\mathsf{vert}}\sfD)$-groupoids.}
\label{fasq}
Following (\ref{defdiag}), a $\sfD$-functor may be viewed as a ($\square^{\mathsf{vert}}\sfD$)-groupoid\,\footnote{One has to check the exactness property.}.

\subsection{Actors.}
\label{act}
It turns out that the basic (algebraic)\,\footnote{Most of these properties own various names in the literature, depending on the authors, and equally unfortunate for our purpose, since these purely algebraic properties received often names issued from Topology, which cannot be kept when working in \textbf{Top} or \emph{Dif}.} properties of a $\sfD$-functor $f:H\rightarrow G$ are encapsulated into \emph{two fundamental squares}, written below, which immediately acquire a ``structured'' (and hence more precise) meaning when written in a diptych.

It is enough to write these properties at the lowest level $(f^{(0)},f^{(1)})$, since it turns out that the special pullback properties of the nerve, as a diagram, allow to carry them over to all levels.

\subsubsection{The activity indicator ${\rm A}(f)$ : in/ex/-actors.}
We consider first the commutative square generated by the \emph{source maps} :

\begin{center}
\begin{diagram}[size=2em,loose,labelstyle=\scriptstyle,inline]
        H     &\rTo^f            &         G   \\
\dSTo^{\alpha_H}         &  {\rm A}(f)              & \dSTo_{\alpha_G}       \\
        E     &\rTo^{f^{(0)}}            &         B   \\
\end{diagram} \ \ .
\end{center}

We shall say $f$ is : an \emph{actor}\,\footnote{See (\ref{GA}) for the terminology.}, an \emph{inactor}, an \emph{exactor}\,\footnote{The underlying algebraic notion is known in the categorical literature under the name of ``fibering functors'' (or ``star-surjective'' functors for Ronnie Brown), which cannot be kept when working in a topological setting.}, depending on whether the square ${\rm A}(f)$ is a pullback, $i$-faithful, or $s$-full (\ref{3bt}).

One can show that any exactor $f:H\rightarrow G$, one can define its \emph{kernel} $K\rITo H$, which is null when $f$ is an actor.

A \emph{principal actor} is an actor $R\stackrel{f}{\rightarrow}G$ with $R$ principal. For instance , taking for $R$ the graph associated to a covering, as described in (\ref{Classat}), we recover the notion of \emph{$G$-cocycle}, including cocycles defining a principal\,\footnote{Whence the terminology.} fibration (when $G$ is a Lie group) and Haefliger cocycles defining a foliation (when G is a pseudogroup).

\subsubsection{The canonical actor $\delta_G$.}
\label{canac}
The map 
$\delta_G:\wedge G=\bigtriangleup G\rSTo G$\ 
may be viewed as a functor, and indeed a principal $s$-actor, associated to the right action of $G$ on itself.

This will be enlightened by the functorial considerations to be developed below.

\subsection{$\sfD$-equivalences.}

\subsubsection{The full/faithfulness indicator ${\rm T}(f)$.}
The second basic square is built with the \emph{transitors} (anchor maps) .
\begin{center}
\begin{diagram}[size=2em,loose,labelstyle=\scriptstyle,inline]
        H     &\rTo^f            &         G   \\
\dTo^{\tau_H}         &  {\rm T}(f)              & \dTo_{\tau_G}       \\
        E\times E     &\rTo^{f^{(0)}\times f^{(0)}}            &         B\times B   \\
\end{diagram} \ \ .
\end{center}
\subsubsection{Equivalences and extensors.}
\label{eqex}
We shall say $f$ is an inductor/$i$-faithful/$s$-full, 
depending on whether the square ${\rm T}(f)$ is a pullback/$i$-faithful/$s$-full\,\footnote{As announced, this explains the terminology used for the squares.} (\ref{3bt}).

When $f$ is $s$-full/a $\sfD$-inductor, and moreover $f^{(0)}$ lies in $\cD_s$, then $f$ is an $s$-functor, and we say $f$ is an \emph{$s$-equivalence}/an \emph{$s$-extensor}
\,\footnote{The terminology derives from the following fact : it turns out that such functors are the exact generalizations of \emph{Lie group extensions} (save for the fact that one has to use two-sided cosets, which, in general, don't coincide with right or left cosets.). The smooth case is treated in [P2], which is written in order to be read possibly in any diptych without any change.}.

While the concept of $\sfD$-inductor derives from the diagrammatic description of ``full and faithful'', the general concept of $\sfD$-equivalence demands to add a diagrammatic description of the ``essential (or generic) surjectivity'', which uses the $A(f)$ square and will not be given here (we refer to [P4]). This general notion allows to speak of $i$-equivalences too.

Parallel to the notion of canonical actor, we have those of canonical equivalences. These will be defined below, when we have given a diagrammatic construction of $\square\,G$ and of the canonical morphisms:
\begin{center}
$G\rightarrow\square\,G\rightrightarrows G\,.$
\end{center}

\subsection{Diptych structures on $\mathbf{Gpd}(\sfD)$.}
\label{diponG}
The following result , which we can just mention here, is of basic importance for a unified study of structured groupoids.

Using the preliminary study of diagrams in a diptych, one can define several useful (Godement) diptych structures on $\mathbf{Gpd}(\sfD)$.

We stress the fact that the $s$-functors are \emph{not} the right candidates for good epis\,\footnote{This problem doesn't arise in the set-theoretic case. Many authors seem to believe that pullbacks along $s$-functors always exist in the $\mathsf{Dif}$ case , but actually the delicate point, often forgotten, is to prove the surmersion condition for the source map.}.

Among various possibilities, one can take:
\begin{itemize}
	\item for \emph{good monos}:
\begin{itemize}
	\item either the $i$-functors
	\item or the $i$-actors;
\end{itemize}
	\item for \emph{good epis} those $s$-functors which moreover belong to one of the following types:
\begin{itemize}
	\item $s$-exactors
	\item $s$-actors
	\item $s$-equivalences.
\end{itemize}
\end{itemize}

\subsection{The category $\mathbf{Gpd}(\sfNc^{+\ast})$.}
\label{cogr}
$\sfD^\ast$-groupoids may be called $\sfD$-cogroupoids.

Some constructions for $\sfD$-groupoids may be better understood from a study of $\sfNc^{+\ast}$-groupoids or $\sfNc^+$-cogroupoids (which are \emph{not} $\mathbf{Set}$-groupoids) (\ref{ficar}).
 Some pieces of notations are needed to avoid confusions arising from duality.

\subsubsection{Notations for $\rm{End}(\cN_c^+)$ and $\rm{End}(\cN_c^{+\ast})$.} 
In contrast to sections (\ref{com}) and (\ref{typeT}), \emph{we shall, for a while, stick to Mac Lane's terminology concerning horizontal and vertical composition of natural transformations} (also called \emph{functorial morphisms}), in order to allow free use of [McL] as reference.

This means that groupoids (viewed as diagrams or functors) have here to be thought as written horizontally (instead of vertically as above), hence the groupoid functors or morphisms (i.e. the arrows of $\mathbf{Gpd}(\sfD)$) as written \emph{vertically}, though this is somewhat uncomfortable.

The natural transformations between endofunctors make up categories denoted by $\mathbf{End}(\cN_c^+)$ and $\mathbf{End}(\cN_c^{+\ast})$, which are indeed double categories (and even more precisely 2-categories [McL]) when considering both horizontal and vertical composition.

The identity maps define canonically: 
\begin{itemize}
	\item a \emph{contravariant} functor:
$$\cN_c^+\rightarrow\cN_c^{+\ast}\,,\,\lambda\mapsto\lambda^\ast\,,$$where $\lambda^\ast$ is $\lambda$ with source and target exchanged ;
 denoting by $\mul$ the product in $\cN_c^{+\ast}$ (i.e. the sum $+$ in $\cN_c^+$), we can write (for any pair of arrows $\lambda$, $\mu$):
$$(\lambda+\mu)^\ast=\lambda^\ast\mul\mu^\ast\,;$$
	\item a \emph{covariant} functor:
$$\rm{End}(\cN_c^+)\rightarrow\rm{End}(\cN_c^{+\ast})\,,\Phi\mapsto{^\ast\Phi^\ast}=\,\stackrel{\ast}{\Phi}\,,$$with $\stackrel{\ast}{\Phi}(\lambda^\ast)=\,(\Phi(\lambda))^\ast$ ;
	\item a bijection
$$\mathbf{End}(\cN_c^+)\rightarrow\mathbf{End}(\cN_c^{+\ast})\,,\,(\varphi:\Phi\stackrel{\cdot}{\rightarrow}\Phi')\mapsto(\varphi^\ast:{\stackrel{\ast}{\Phi}}\,\stackrel{\cdot}{\leftarrow}\,{\stackrel{\ast}{\Phi}}{'})\,,$$with $\varphi^\ast(n)=(\varphi(n))^\ast$, which is: 
	\begin{itemize}
	\item \emph{covariant} with respect to \emph{horizontal} composition 
laws ;
	\item \emph{contravariant} with respect to \emph{vertical} composition 
laws.
	\end{itemize}
\end{itemize}

\subsubsection{Description of $\mathbf{Gpd}(\sfNc^{+\ast})$.}
\label{exng}
By the previous bijection:
\begin{itemize}
	\item $\sfNc^{+\ast}$-groupoids derive \emph{covariantly} from the endofunctors of $\sfNc^+$ \emph{preserving surjections, injections and pushouts} ;
	\item morphisms of $\sfNc^{+\ast}$-groupoids derive \emph{contravariantly} from morphisms between endofunctors of the previous type.
\end{itemize}

Moreover we know from the general descriptions of $\sfD$-groupoids (\ref{codes}) that:
\begin{itemize}
	\item the endofunctors $\Gamma={\stackrel{\ast}{\Phi}}$ of $\sfNc^{+\ast}$ defining $\sfNc^{+\ast}$-groupoids are \emph{uniquely} determined by the data:$$(\Gamma^{(0)},\Gamma^{(1)},\omega_\Gamma,\alpha_\Gamma,\delta_\Gamma),$$
	hence (resulting from the previous study) by the data in $\sfNc^+$:
	$$(n_0=\Phi(\cdot0),n_1=\Phi(\cdot1),n_2=\Phi(\cdot2),\omega:n_0\lSTo n_1,\alpha:n_1\lITo n_0,\delta:n_2\lITo n_1).$$
\end{itemize}

\subsubsection{Some basic examples of $\sfNc^{+\ast}$-groupoids and morphisms.}
\label{dig}
As just explained, such a groupoid morphism $\gamma:\Gamma\leftarrow\Gamma'$ derives from a functorial morpism $\varphi:\Phi\stackrel{\cdot}{\rightarrow}{\Phi'}$ where $\Phi,\Phi'$ preserve surjections, injections and pushouts. These define a 2-subcategory of $\mathbf{End}(\cN_c^+)$ denoted by $\digamma$.

Those $\Phi$\,'s which preserve sums too are of type :
\begin{center}
$\boldsymbol{p_\times}:\cdot n\mapsto p\times\cdot n,\,\lambda\mapsto p\times \lambda$ ($p$ times iterated sum in $\sfNc^+$).
\end{center}
The associated groupoids, denoted by $\boldsymbol{p_\mul}$ are the \emph{banal} (but not principal\,\footnote{Since the terminal object 0 has been dropped.}) $\sfNc^{+\ast}$-groupoids, in the sense of (\ref{bagr})\,\footnote{While the banal (and principal) $\sfNc^+$-groupoids are $p\times p$.}\,.

For $p=0$, we get the ``\emph{null groupoid}'', denoted by $\mathbf0$, associated to the constant functor $\mathbf0:\cdot n\mapsto\cdot0$, and, for $p=1$, the\, ``\emph{unit groupoid}'', denoted by $\mathbf1$, associated to the identity functor.

For $p=2$, one has the ``\emph{square groupoid}''\,\footnote{The terminology will become clear below in \ref{exr}.}\,, denoted by $\boldsymbol{\cast}=\mathbf2_\mul$, associated to:
$$\boldsymbol{2_\times}:\cdot n\mapsto2(\cdot n)=
\cdot n+\cdot n,\,\lambda\mapsto2\lambda=\lambda+\lambda\,.$$

Among the non principal ones (\ref{palg}), the simplest is denoted by $\boldsymbol\tast$ : it is associated to the ``\emph{shift functor}'':$$\boldsymbol{\cdot0_+}:\cdot n\mapsto\cdot0+\cdot n,\,\lambda\mapsto\cdot0+\lambda$$\,
where the last $\cdot0$ is understood as the identity of the object $\cdot0=1$.

We have also a $\sfNc^{+\ast}$-groupoid morphism: 
\begin{center}(1)\quad
\fbox{
$\boldsymbol{\delta_0^\ast:\,\tast\longrightarrow1}$}\,
\end{center}
associated, with the notations of (\ref{groudi}) and of [McL], to the the ``\emph{shift morphism}'', defined by the family $(\delta_0^{\cdot n})_{n\in\bbN}$ (injections skipping the $0$ in $\cdot1+\cdot n$).

As to $\boldsymbol\cast$\,, associated to $\boldsymbol{2_\times}$, we have in $\sfNc^+$ the coproduct morphisms $(n\in\bbN)$ :
\begin{center}
{\begin{diagram}[size=3.5em,tight,labelstyle=\scriptstyle,inline]
(\cdot n)&\lTo^{\rm{codiag}}&(\cdot n)+(\cdot n)&\pile{\lTo^{\quad\iota_2}\\ \lTo_{\quad\iota_1}}&(\cdot n)\\
\end{diagram}}
\end{center}
which define morphisms of $\sfNc^{+\ast}$-groupoids denoted suggestively by:

\begin{center}(2)\quad
\fbox{
\begin{diagram}[size=2.2em,tight,labelstyle=\scriptstyle]
\boldsymbol{1}&\rTo^{\boldsymbol{\omega}}&\boldsymbol{\cast}&\pile{\rTo^{\boldsymbol{\varpi_2}}\\ \rTo_{\boldsymbol{\varpi_1}}}&\boldsymbol{1}\\
\end{diagram}}\,.
\end{center}


By iteration we can even get a \emph{canonical groupoid $(\cast^{(n)})_{n\in\bbN}$ in $\mathbf{Gpd}(\sfNc^{+\ast})$}\,\footnote{More precisely, with a suitable diptych structure (see \ref{diponG}). This is indeed a double groupoid, or better a ``\emph{groupoid-cogroupoid}''.}.
\par\smallskip
The ``\emph{symmetry groupoid}'' $\boldsymbol\Sigma_\ast$ proceeds from the (involutive) ``\emph{symmetry map}'' $\Sigma\,:\cN_c^+\rightarrow\cN_c^+$, derived by reversing the order on the integers viewed as ordinals $(\overrightarrow n\mapsto\overleftarrow n)$ ; it is defined by the identity on the objects, and, on the generators (with Mac Lane's notations), by $(n\in\bbN^+)$:
\begin{center}
\fbox{$\Sigma:\delta^n_j\mapsto\delta^n_{n-j},\,\sigma^n_j\mapsto\sigma^n_{n-j}$}\,,
\end{center}
and the family of maps : $(\varsigma_n:n\rightarrow n,\,j\mapsto n-j)_{(n\in\bbN^+)}$ defines a natural transformation from identity towards $\Sigma$, hence also a groupoid morphism:
\begin{center}(3)\quad
\fbox{$\boldsymbol{\varsigma:\Sigma_\ast\longrightarrow1}$\,.}
\end{center}

\section{Double functoriality of the definition of $\sfD$-groupoids.}

\subsection{Bivariance of $\sfD$-groupoids and $\sfD$-functors.}
The definition of a $\sfD$-groupoid as a
\begin{diagram}[size=4em,labelstyle=\scriptscriptstyle,textflow]
\sfNc^{+\ast}&
\pile{\rTo^{\Gamma}\\
 \scriptscriptstyle{\downarrow\gamma}\\
 \rTo_{\Gamma'}}&
\sfNc^{+\ast}&
\pile{\rTo^{\mathbf H}\\
 \scriptscriptstyle{\downarrow\mathbf f}\\
 \rTo_{\mathbf G}}&
\sfD&
\pile{\rTo^{\mathsf T}\\
 \scriptscriptstyle{\downarrow\mathsf t}\\
 \rTo_{\mathsf{T'}}}
&\sfD'\\
\end{diagram}
 diptych morphism $\mathbf G :\sfNc^{+\ast}\rightarrow \sfD$ shows immediately (look at the adjoining diagram, where we go on sticking to Mac Lane's conventions) this definition is:
\begin{itemize}
	\item functorial with respect to composition, on the left (target side), with an exact diptych morphism $\sfT:\sfD\rightarrow\sfD'$\,, as well as with respect to right (horizontal ) composition with a natural transformation $\mathsf{t:T\stackrel{\cdot}{\rightarrow}T'}$;this means that
any $\sfD$-functor $\mathbf{f:H\rightarrow G}$ gives rise to a $\sfD'$-functor $$\mathsf{T}\circ\mathbf{f}:\mathsf{T}\circ\mathbf{H}\rightarrow \mathsf{T}\circ\mathbf{G}$$\footnote{As in [McL], the same notation is used for a functor and the identity natural transformation associated to this functor, here $\mathsf{T}$} and to a natural transformation:
$$\mathsf{t}\circ\mathbf{f}:\mathsf{T}\circ\mathbf{f}\stackrel{\cdot}{\rightarrow}\mathsf{T'}\circ\mathbf{f}\,;$$
the latter gives rise to the following commutative square of $\sfD'$-groupoids\,\footnote{Written in loose notations, i.e. identifying groupoids and functors with their 1-level part.}:
\begin{center}
\begin{diagram}[size=2.5em,labelstyle=\scriptstyle,inline]
TH&\rTo^{Tf}&TG\\
\dTo^{t(H)}&&\dTo_{t(G)}\\
T'H&\rTo^{T'f}&T'G\\
\end{diagram}
\ \ ;
\end{center}
actually $\mathbf{Gpd}$ behaves here like a functor and we can define:
\begin{center}
$\mathbf{Gpd}(\mathsf T):\mathbf{Gpd}(\sfD)\rightarrow\mathbf{Gpd}(\sfD')$, 
$\mathbf{Gpd}(\mathsf t):
\mathbf{Gpd}(\mathsf T)\stackrel{\cdot}{\rightarrow}\mathbf{Gpd}(\mathsf T')$\,;
\end{center}

	\item functorial with respect to horizontal composition on the right (source side), with a $\sfNc^{+\ast}$-groupoid $\Gamma$, as well as with a groupoid morphism $\gamma\,:\,
\Gamma\rightarrow\Gamma'$\, viewed as a natural transformation between endofunctors of $\sfNc^{+\ast}$\,; this means that any $\sfD$-functor $\mathbf{f:H\rightarrow G}$ gives rise to a $\sfD$-functor: 
\begin{center}
$\mathbf{\Gamma^\bullet f\,:
\Gamma^\bullet H\rightarrow\Gamma^\bullet G}$,
\end{center}
where 
\begin{center}
$\mathbf{\Gamma^\bullet H=H\circ\Gamma}$, 
$\mathbf{\Gamma^\bullet G=G\circ\Gamma}$, 
$\mathbf{\Gamma^\bullet f=f\circ\Gamma}$\,,
\end{center}
and to a natural transformation
\begin{center}
$\mathbf{\gamma^\bullet f=f\circ\gamma\,:
\Gamma^\bullet G\rightarrow\Gamma'^\bullet G}$ ;
\end{center}
this defines a (covariant) functor:
$$\mathbf{\Gamma^{\bullet}:Gpd(\sfD)\rightarrow Gpd(\sfD)}\,,$$
hence a \emph{canonical \emph{(vertical)} representation}:
\begin{center}
\fbox{
$\mathbf{Gpd}(\sfNc^{+\ast})\rightarrow\mathbf{Gpd}(\sfD)$}\,,
\end{center}
but this representation depends in a \emph{contravariant} way upon $\Gamma$, with respect to the \emph{horizontal} composition of $\sfNc^{+\ast}$-groupoids defined above, since 
\begin{center}
$\mathbf{\Gamma^{\bullet}\circ\Gamma'^{\bullet}=(\Gamma'\circ\Gamma)^{\bullet}}$\,,\ 
$\mathbf{\gamma^{\bullet}\circ\gamma'^{\bullet}=(\gamma'\circ\gamma)^{\bullet}}$\,.
\end{center}
\end{itemize}
\par\smallskip
Note that, when going back to the generating endofunctors of $\sfNc^+$, we get a \emph{doubly contravariant canonical representation} of $\digamma$ (notation of \ref {dig}) into $\mathbf{Gpd}(\sfD)$.

\par\smallskip
We give a few examples.

\subsection{Examples for the left functoriality.}
We can either ``forget'' the structures on objects of $\sfD$, or ``enrich'' them. We give examples of these opposite directions.

\subsubsection{Concrete diptychs.}

Thinking to the basic examples of $\mathsf{Top}$ and $\mathsf{Dif}$, we shall say $\sfD$ is ``\emph{concrete}'' if it comes equipped with an adjunction from $\sfE=\mathsf{Set}$ to $\sfD$ defined by the following adjoint pair [McL] of functors\,\footnote{Assumed moreover to be faithful, to preserve products, and to define exact diptych morphisms. An object $B$ of $\cD$ is viewed as a ``structure'' on the ``underlying set'' $|B|$, and an arrow $A\rightarrow B$ of $\cD$ is fully described by the triple $(A,|f|,B)$\,. Any set $E$ may be endowed with the ``discrete structure'' $\dot{E}$.}:
\begin{center}
(discrete, forgetful) $=(\ {\dot{}}\ ,|\,|)$\ :\ \begin{diagram}[size=1.8em,labelstyle=\scriptscriptstyle,inline]
\sfE&\pile{\rTo^{\dot{}}\\ \lTo_{|\,|}}&\sfD\\
\end{diagram}
\,.
\end{center}

Then we can speak of the underlying $\mathsf{E}$-groupoid, and make use of set-theoretical descriptions.

\subsubsection{The tangent functor.}

Thinking now to the case when $\sfD=\mathsf{Dif}$, we can consider (see \ref{vb}) the tangent functor:
$$\mathsf{T:Dif\longrightarrow VectB}\,,$$
which is equipped with two natural transformations:
$$0\rTo^\cdot_o T\rTo^\cdot_t 0,.$$

Once one has checked it defines an exact diptych morphism, we can immediately transfer to the tangent groupoids all the general constructions valid for general $\sfD$-groupoids (for instance constructions of fibred products, and so on).

\subsubsection{Double groupoids.} The idea of defining and studying the notion of double groupoids as groupoids in the category of groupoids is due to Ehresmann, who proved the equivalence with the alternative description by means of two category composition laws satisfying the ``exchange law''\,\footnote{In our framework this would result from \ref{exp}.}.

In the diptych setting we get several notions depending on the choice for the diptych structure on $\rm{Gpd}(\sfD)$ (see \ref{diponG}).

We stress the point that, even in the purely set-theoretical setting, the choice we made of exactors for good epis implies adding a certain surjectivity condition\,\footnote{If we consider the square made up by the sources and targets of both laws, this condition means that three of the four edges may be given arbitrarily.} which does not appear in Ehresmann's definition, but was encountered by several authors, mainly Ronnie Brown (filling condition), and seems useful to develop the theory beyond just definitions.

Applying the general theory of $\sfD$-diptychs, for instance for getting fibred products, or quotients, or Morita equivalences, gives results which it would be very hard to get by a direct study (which has never been done), even in the purely algebraic setting.

\subsection{Examples for the right functoriality.}
\label{exr}
As announced,the coherence of various notations introduced above will appear just below (to be precise, we find : ${\tast\,}^\bullet=\bigtriangleup,\,{\cast}\:{^\bullet}=\square$).
The examples given in \ref{dig}, allow to transfer to $\sfD$-groupoids some classical set-theoretic constructions, announced above. Taking for $\gamma:\Gamma\rightarrow\Gamma'$:
\begin{itemize}
	\item formula (1) of \ref{dig}, we recover the morphism: 
	$$\delta_G:\bigtriangleup G\longrightarrow G\,;$$
	\item formula (2), we get the canonical equivalences:
	$$G\rightarrow\square\,G\rightrightarrows G\,;$$
	\item formula (3), we get the dual groupoid $G^\ast=\Sigma G$, and the \emph{inverse law} $\varsigma_G=\iota_G$, which defines an involutive isomorphism:
	$$\varsigma_G:G\rightarrow G^\ast\,.$$
\end{itemize}
\par\smallskip

\subsection{$\sfD-$natural transformations, holomorphisms.}
\label{hol}
Once $\square G$ is defined, one can also define $\sfD$-\emph{natural transformations} between functors from $H$ to $G$ as $\sfD$-functors $H\rightarrow\square G$. The canonical $i$-equivalence $G\rightarrow\square G$ defines the identical transformation of the identity functor.

Since $G$ is a groupoid, such natural transformations are necessarily functorial isomorphisms. An isomorphy class of $\sfD$-functors will be called a ``\emph{holomorphism}'' (an alternative terminology might be ``exomorphism'', since this notion generalizes the outer automorphisms of groups).

Since the horizontal composition of natural transformations commutes with the vertical one, it defines a composition between holomorphisms, and this yields a quotient category of $\textbf{Gpd}(\sfD)$, which we shall denote by $\textbf{Hol}(\sfD)$.

Any $\sfD$-functor $f:H\rightarrow G$ generates the following commutative diagram,
\begin{diagram}[size=1.8em,tight,textflow,midshaft,labelstyle=\scriptstyle,midshaft]
&&G\\
&\ruTo^{p(f)}&\uSTo^{\pi_2}_\sim\\
K&\rTo&\square G\\
\dSTo^{q(f)}_\sim&\pb&\dSTo^{\pi_1}_\sim\\
H&\rTo^f&G\\
\end{diagram}
 in which $q(f)$ is an $s$-equivalence, and $p(f)$ an exactor. The pair $(p(f),q(f))$ is called the ``\emph{holograph}'' of $f$.

Moreover the $s$-equivalence $q(f)$ is \emph{split} (\ref{Dfun}) since $\pi_1$ is.

\section{The butterfly diagram.} This section illustrates, in the case of orbital structures presently described, the use of diagrams in a diptych for transferring constructions in \textbf{Set} to constructions in \textbf{Dif} or other various categories, as well as the use of the various kinds of $\sfD$-functors introduced above. We can be only very sketchy. More precise descriptions and results may be found in [P4], where they are stated for the differentiable case, but written to be easily transferable to general diptychs. More details about the diagrams used for proofs will be given elsewhere.
\subsection{Generalized ``structure'' of the orbit space.}
\subsubsection{Algebraic ``structure''.}
First, from a purely set-theoretic point of view, any equivalence between two groupoids, in the general categorical sense [McL]\,\footnote{i.e. a functor which is full, faithful, and essentially surjective, but possibly non-surjective.}, preserves\,\footnote{More precisely, this means that it defines a bijection between the two orbit spaces.} the set-theoretic orbit space $Q$, but indeed it preserves much more, since, to each orbit (or transitive component), there is a (well defined only up to isomorphism) attached isotropy group, and this endows $Q$ with a kind of an algebraic ``structure'', in a (non set-theoretical) generalized sense.
For instance, the leaves of a foliation are marked by their holonomy groups, the orbits of a group(oid) action are marked by their fixators. Such a ``structure'' is sometimes called a group stack.

\subsubsection{$\sfD$\,-``structure'' of the orbit space.}
Then, in the $\sfD$-framework, replacing algebraic equivalences by $\sfD$-equivalences will moreover encapsulate in this generalized structure the memory of the $\sfD$-structure as well. It turns out (though this is by no means \emph{a priori} obvious) that it is enough to make use of $s$-equivalences (\ref{eqex}). Finally we are led to the following:

\textbf{Definition:} Two $\sfD$-groupoids $H$, $G$, are said to be \emph{$\sfD$-equivalent} if they are linked by a pair of $s$-equivalences:
\begin{diagram}[size=2em,tight,midshaft,inline,labelstyle=\scriptstyle]
H&\lSTo^q_\sim&K&\rSTo^p_\sim&G
\end{diagram}\,.

The fact that this is indeed an equivalence relation is an easy consequence of the results stated in \ref{diponG}, taking $s$-equivalences as good epis in $\mathbf{Gpd}(\sfD)$, and using fibred products of good epis.

A $\sfD$-equivalence class of $\sfD$-groupoids may be called an \emph{orbital structure}. Any representative of this equivalence class is called an \emph{atlas} of the orbital structure. 

\subsection{Inverting equivalences.}
\subsubsection{Meromorphisms.}
Note than in general orbital structures cannot be taken as \emph{objects} of a new category. However one can define [P4] a new category which shall be denoted here by $\textbf{Mero}(\sfD)$, with the \emph{same objects} as $\textbf{Gpd}(\sfD)$, and arrows called ``\emph{meromorphisms}'', in which the $s$-equivalences (and indeed all the $\sfD$-equivalences) become invertible, in other words are turned into isomorphisms. In the topological case, these isomorphisms may be identified with the \emph{Morita equivalences}. In this new category the orbital structures now become \emph{isomorphy classes} of $\sfD$-groupoids, though the objects still remain $\sfD$-groupoids and not isomorphy classes, so that the orbital ``structures'' \emph{are not carried by actual sets}, and, as such, remain ``virtual''.

This means that $\textbf{Mero}(\sfD)$, is the universal solution for the problem of fractions consisting in formally inverting the $s$-equivalences, and indeed all the $\sfD$-equivalences.

This kind of problem always admits a general solution [G-Z]: the arrows are given by equivalence classes of diagrams, the description of which is simpler when the conditions for right calculus of fractions are satisfied.

It is worth noticing that these assumptions are not fulfilled here, while our construction is in a certain sense much simpler, since, as is the case in Arithmetic, each fraction will admit of a simplified or irreducible canonical representative $(p,q)$. In the topological case these representatives may be identified with the ``generalized homomorphisms'' described by A. Haefliger in [Ast 116] and attributed to G. Skandalis, or the ``$K$-oriented morphisms'' of M. Hilsum and G. Skandalis.

We can say more. The canonical functor $\textbf{Gpd}(\sfD)\stackrel{\Phi}{\rightarrow}\textbf{Mero}(\sfD)$
admits of the following factorization:
$\textbf{Gpd}(\sfD)\stackrel{\Phi_1}{\rightarrow}\textbf{Hol}(\sfD)\stackrel{\Phi2}{\rightarrow}\textbf{Mero}(\sfD)$
(see \ref{hol}) with $\Phi_1$ full (i.e. here surjective) and $\Phi_2$ faithful (or injective).

It turns out that $\textbf{Hol}(\sfD)$ is the solution of the problem of fractions for \emph{split} (\ref{Dfun}) $\sfD$-equivalences. It is embedded in $\textbf{Mero}(\sfD)$ by means of the holograph (\ref{hol}).

\subsubsection{Description of fractions.}
Let $(p,q)$ denote a pair of \emph{exactors} with the same source $K$:
$p:K\rightarrow G,\,q:K\rightarrow H$. We set $R=\text{Ker}\,q,S=\text{Ker}\,p$ (see \ref{act}).

Let $(p',q')$ another pair $p':K'\rightarrow G,\,q':K'\rightarrow H$ with the same $G$ and $H$.

Letting for a while $G$ and $H$ fixed, we start considering arrows $k:(p',q')\rightarrow(p,q)$ defined as $D$-functors $k:K'\rightarrow K$ such that the whole diagram commutes.

Let $p/q$ denote the isomorphy class of $(p,q)$, and call it a ``\emph{fraction}''.

On the other hand we say two pairs $(p_i,q_i)(i=1,2)$ are \emph{equivalent} if there exist two $s$-equivalences $k_i:(p,q)\rightarrow(p_i,q_i)$. This is indeed an equivalence relation, and the class of $(p,q)$ will be denoted by $pq^{-1}$.

We consider now those pairs $(p,q)$ satisfying the subsequent extra conditions, which turn out to be preserved by the previous equivalence:
\begin{enumerate}
	\item $q$ is an $s$-equivalence;
	\item $p$ and $q$ are ``\emph{cotransversal}''.
\end{enumerate}
The former condition implies that the kernel $R=\text{Ker}\,q$ is principal (\ref{palg}). 
The latter condition will be expressed by means of the following (commutative) ``\emph{butterfly diagram}'' gathering the previous data:

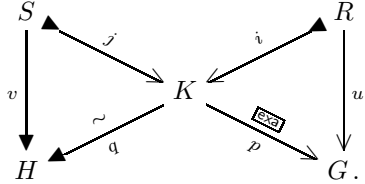
\begin{diagram}[w= 3em,h=1.5em,tight,labelstyle=\scriptstyle]
S&&&&R\\
&\rdITo^j&&\ldITo^i&\\
\dSTo^v&&K&&\dTo_u\\
&\ldSTo^\sim_q&&\rdTo^\exa_p&\\
H&&&&G\,.\\
\end{diagram}
Then the condition of cotransversality means that $u$ or (this is indeed equivalent) $v$ is an exactor (then $v$ will be an $s$-exactor).

When $u$ (and $v$) are actors, $p$ and $q$ are said to be ``\emph{transverse}'', and the fraction $p/q$ is called ``\emph{irreducible}'' (or simplified).

One can show (using the theory of extensors) that the class $pq^{-1}$ owns a \emph{unique irreducible representative} $p/q$.

Our \emph{meromorphisms} (from $H$ to $G$) are then defined as the classes $pq^{-1}$ or their irreducible representatives $p/q$.

A \emph{``Morita equivalence''} is the special case when $p$ is an $s$-equivalence too. The butterfly diagram is then \emph{perfectly symmetric}. We say that $(u,v)$ is a pair of ``\emph{conjugate principal actors}''. Each one determines the other one up to isomorphism.

Using irreducible representatives and forgetting the $\sfD$-structures, one then recovers easily the set-theoretical part of the description of Skandalis-Haefliger homomorphisms (two commuting actions, one being principal). Now, in the differentiable case, the local triviality conditions are automatically encapsulated in the surmersion conditions (imposed to the good epis) by means of the Godement theorem.

One of the immense advantages of this presentation (apart from being defined in many various frameworks), is that the use of non irreducible representatives allows a very natural definition of the composition of meromorphisms (note that in [Ast 116] this composition is defined by A. Haefliger but in very special cases, when one arrow is a Morita equivalence). This composition is defined by means of the following diagram (using the diptych properties of $\textbf{Gpd}(\sfD)$) (\ref{diponG}):
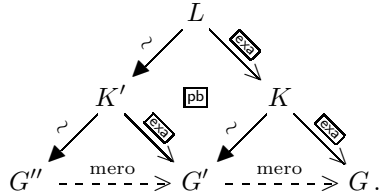
\begin{diagram}[size=1.6em,tight,labelstyle=\scriptstyle,midshaft]
&&&&L&&&&\\
&&&\ldSTo^{\sim}&&\rdTo^\exa&&&\\
&&K'&&\pb&&K&&\\
&\ldSTo^{\sim}&&\rdTo^\exa&&\ldSTo^{\sim}&&\rdTo^\exa&\\
G''&&\rDashTo^{\text{mero}}&&G'&&\rDashTo^{\text{mero}}&&G\,.\\
\end{diagram}
(Of course there are many things to check to justify all our claims).

\subsection{Example.}
Let us come back to the example of the space of leaves of a (regular) foliation (\ref{Fol}). We invite the reader to look at what happens when we take as good epis:
\begin{enumerate}
	\item all the surmersions;
	\item the retroconnected ones;
	\item the retrodiscrete (or \emph{étale}) ones;
	\item the proper ones.
\end{enumerate}
In the first case, we are allowed to take as an atlas the holonomy groupoid and the transverse holonomy pseudogroups associated to various tranversals as well, which all belong to the same Morita class.

The second choice is adapted to the search for invariants of the Molino equivalence class of a foliation: the Morita class of the holonomy groupoid is such an invariant.

The third choice is adapted to the use of holonomy pseudogroups and of van Est S-atlases [Ast 116], and to the study of the effect of coverings on foliations.

The fourth choice is adapted to the study of compact leaves of foliations and of properties related to Reeb stability theorem, as well as to the study of orbifolds (or Satake manifolds).

In fact it is useful to use various choices simultaneously.

\section{Epilogue.}
\begin{itemize}
	\item We have not given here statements concerning diptych structures on the category $\textbf{Mero}(\sfD)$, since our present results are still partial and demand some further checks to ensure them completely, especially concerning the Godement property. It is clear for us that such types of statements would be very useful, since, for instance, groupoids in such diptychs would be fascinating objects. Anyway it seems clear for us that this category has to be explored more deeply.
	\item We know that our formal construction for the previous category of fractions seems to work perfectly as well when replacing the category of $s$-equivalences by that of $s$-extensors or by various subcategories of the latter (adding for instance conditions of connectedness on the isotropy groups). We are convinced that such categories, which are much less known (not to say totally unknown) than the previous one, are basic for the understanding of holonomy of foliations with singularities ($\check{S}tefan\,foliations$), and that these enlarged Morita classes certainly encapsulate some deep and hidden properties of the orbit spaces.
\end{itemize}

\small\textbf{Acknowledgements.}
\small\begin{enumerate}
	\item I would like to thank warmly the organizers for their invitation to this Conference at Krynica. I was happy with the friendly and stimulating atmosphere which was reigning throughout this session.
	\item I am indebted to Paul Taylor for the (certainly very awkward) use I made of his package ``diagrams'', allowing various styles of arrows, which are very useful as condensed visual mnemonics for memorizing properties of maps and functors.
\end{enumerate}


\bibliographystyle{amsplain}
\providecommand{\bysame}{\leavevmode\hbox to3em{\hrulefill}\thinspace}

\end{document}